\documentclass[12pt]{article}                
\usepackage{graphicx}
\usepackage{psfrag}
\usepackage{amsmath, amssymb}
\usepackage{mathptmx}
\usepackage{ntheorem}
\usepackage[nottoc,numbib]{tocbibind}
\usepackage{hyperref}

\usepackage{setspace}
\makeatletter
\newsavebox\myboxA
\newsavebox\myboxB
\newlength\mylenA
\newcommand*\xoverline[2][0.65]{%
    \sbox{\myboxA}{$\m@th#2$}%
    \setbox\myboxB\null
    \ht\myboxB=\ht\myboxA%
    \dp\myboxB=\dp\myboxA%
    \wd\myboxB=#1\wd\myboxA
    \sbox\myboxB{$\m@th\overline{\copy\myboxB}$}
    \setlength\mylenA{\the\wd\myboxA}
    \addtolength\mylenA{-\the\wd\myboxB}%
    \ifdim\wd\myboxB<\wd\myboxA%
       \rlap{\hskip 0.6\mylenA\usebox\myboxB}{\usebox\myboxA}%
    \else
        \hskip -0.6\mylenA\rlap{\usebox\myboxA}{\hskip 0.6\mylenA\usebox\myboxB}%
    \fi}
\makeatother

\newcommand{\A}{{\cal A}}

\newcommand{\DD}{{\cal D}}
\newcommand{\PP}{{\cal P}}

\newcommand{\tfi}{{\tilde{\Phi}}}

\newcommand{\eps}{{\varepsilon}}
\newcommand{\ch}{{\mbox{\rm ch}}}

\newcommand{\smsp}{\hspace{0.3mm}}
\newcommand{\e}{\mathbb{E}}

\newcommand{\p}{\mathbb{P}}

\newcommand{\Reals}{\mathbb{R}}

\newcommand{\la}{\langle}
\newcommand{\ra}{\rangle}
\newcommand\qed{\hfill\hbox{\rlap{$\sqcap$}$\sqcup$}}

\newtheorem{proposition}{Proposition}
\newtheorem{lemma}{Lemma}
\newtheorem{theorem}{Theorem}

\theoremstyle{nonumberplain}
\newcommand\specialref{}

\begin{document}

\title{Chaos in temperature in generic $2p$-spin models}
\author{Dmitry Panchenko\thanks{\textsc{\tiny Department of Mathematics, University of Toronto, panchenk@math.toronto.edu. Partially supported by NSERC.}}
}
\date{}
\maketitle
\begin{abstract}
We prove chaos in temperature for even $p$-spin models which include sufficiently many $p$-spin interaction terms. Our approach is based on a new invariance property for coupled asymptotic Gibbs measures, similar in spirit to the invariance property that appeared in the proof of ultrametricity in \cite{PUltra}, used in combination with Talagrand's analogue of Guerra's replica symmetry breaking bound for coupled systems.
\end{abstract} 
\vspace{0.5cm}

\section{Introduction}

The phenomenon of chaos in temperature in spin glass models was first studied in the physics literature by Fisher and Huse \cite{FH86} and Bray and Moore \cite{BM87} and can be briefly described as follows. One can show in some spin glass models, such as the Sherrington-Kirkpatrick and mixed $p$-spin models, that the Gibbs distribution of the system at a given temperature is concentrated near some constant level of energy (at the right scale) and this level can change with temperature. Chaos in temperature means that the set of likely configurations looks quite different even if we change temperature only slightly, and if we sample two spin configurations from the Gibbs distributions at different temperatures then their overlap will be almost deterministic. It means that the two systems might have a common preferred direction, for example in the presence of external field, but otherwise are completely uncorrelated. This is in contrast with the behaviour of the overlap of two configurations from the system at the same temperature, which may have a non-trivial distribution according to the Parisi ansatz (see \cite{MPV, SKmodel, SG2-2}). In this paper, we will prove chaos in temperature for mixed even $p$-spin models that include sufficiently many $p$-spin interactions.

We will consider a mixed $p$-spin model, which is a generalization of the Sherrington-Kirkpatrick model \cite{SK}, corresponding to the Hamiltonian
\begin{equation}
H_N(\sigma) = \sum_{p\geq 2} \gamma_p H_{N,p}(\sigma)
\label{Hammixed}
\end{equation}
defined on $\{-1,+1\}^N$, where the $p^{\mathrm{th}}$ term
\begin{equation}
H_{N,p}(\sigma)
=
\frac{1}{N^{(p-1)/2}}
\sum_{i_1,\ldots,i_p = 1}^N g_{i_1\ldots i_p} \sigma_{i_1}\cdots\sigma_{i_p}
\label{ch11mixedp}
\end{equation}
is called a pure $p$-spin Hamiltonian, coefficients $(g_{i_1\ldots i_p})$ are standard Gaussian random variables independent for all $p\geq 2$ and all $(i_1,\ldots,i_p)$, and coefficients $(\gamma_p)_{p\geq 2}$ decrease fast enough, for example, $\sum_{p\geq 2} 2^p \gamma_p^2<\infty$, to ensure that the Hamiltonian is well defined when the sum includes infinitely many terms. An important feature of these models is that, if we denote by
\begin{equation}
R_{1,2} = \frac{1}{N} \sum_{i=1}^N \sigma_i^1  \sigma_i^2 
\label{ch11overlap}
\end{equation}
the overlap of two configurations $\sigma^1,\sigma^2\in \{-1,+1\}^N$, then the covariance of the Gaussian process $H_N(\sigma)$ in (\ref{Hammixed}) is a function of the overlap,
\begin{equation}
\e H_N(\sigma^1) H_N(\sigma^2) = N\xi(R_{1,2}), 
\label{ch11xidefine}
\end{equation}
where 
$
\xi(x)=\sum_{p\geq 2}\gamma_p^2 x^p.
$
In this article we will only consider \emph{generic} even $p$-spin models defined as follows.

\medskip
\noindent
\textbf{Definition.} \emph{(Generic even $p$-spin model)} We will call the above mixed $p$-spin Hamiltonian \emph{generic} if $\gamma_p=0$ for odd $p\geq 3$ and linear span of functions $x^p$ for even $p\geq 2$ such that $\gamma_p\not = 0$ is dense in $C([0,1],\|\,\cdot\,\|_\infty)$.

\medskip
\noindent
In other words, we assume that sufficiently many $p$-spin interaction terms are included in the Hamiltonian of the model. By the M\"untz-Sz\'asz theorem, the density condition is equivalent to $\sum_{p\geq 1} p^{-1}I(\gamma_p\not = 0)=\infty.$ This will be needed to obtain the general form of the Ghirlanda-Guerra identities \cite{GuerraGG, GG} that will be used crucially in the proof of our main result. 

When we consider two copies of the system, we can include arbitrary external fields, so let us consider a random vector $(h^1,h^2)$ and i.i.d. copies $(h_i^1,h_i^2)_{i\geq 1}$. The distribution of $(h^1,h^2)$ can be arbitrary, and not necessarily centered. We only need some integrability condition and will assume that both coordinates have subgaussian tails. Then we consider the Hamiltonians
\begin{equation}
H_N^j(\sigma) = H_N(\sigma)+ \sum_{i\leq N}h_i^j\sigma_i
\label{plush}
\end{equation}
for $j=1,2.$
Given inverse temperature parameters $\beta_1>0$ and $\beta_2>0$,
\begin{equation}
G_{N}^j(\sigma) = \frac{\exp \beta_j H_N^j(\sigma)}{Z_{N}^j},
\,\,\mbox{ where }\,\,
Z_{N}^j = \sum_{\sigma} \exp \beta_j  H_N^j(\sigma),
\label{ch12G}
\end{equation}
denotes the Gibbs measure of the $j^{\mathrm{th}}$ system, where $Z_{N}^j$ is called the partition function. 

We will denote by $(\sigma,\rho)$ the vectors in $(\{-1,+1\}^N)^2$ and by $\la\,\cdot\,\ra$ the average with respect to $(G_{N}^1\times G_{N}^2)^{\otimes\infty}$. We will use the notation $\tilde{\sigma}=\sigma/\sqrt{N}$ and $\tilde{\rho}=\rho/\sqrt{N}$, so that the overlap in (\ref{ch11overlap}) can be written as $\tilde{\sigma}^1\cdot\tilde{\sigma}^2$. The overlap between a replica $\sigma^1$ from $G_N^1$ and replica $\rho^1$ from $G_N^2$ can be written as $\tilde{\sigma}^1\cdot \tilde{\rho}^1$. Our main result is the following.

\begin{theorem}\label{Th1}
If $\beta_1\not = \beta_2$ then there exists a constant $\chi\in[-1,1]$ such that
\begin{equation}
\lim_{N\to\infty}\e \bigl\la \bigl(\tilde{\sigma}^1\cdot \tilde{\rho}^1- \chi \bigr)^2\bigr\ra = 0.
\label{Th1Eq}
\end{equation}
Moreover, $\chi=0$ if either $\e(h^1)^2=0$ or $\e(h^2)^2=0$.
\end{theorem}
A result of this nature was proved by Chen in \cite{ChenChaos}, but it required tuning the parameters $(\gamma_p)$ in (\ref{Hammixed}) between the two systems in a special way instead of changing the global inverse temperature parameter $\beta$ as in Theorem \ref{Th1}, which is the canonical form of chaos in temperature. Previous results concerned with another type of chaos, disorder chaos, were obtained by Chatterjee in \cite{Chatt09} (see \cite{Chatt14}) in the case of no external field, and by Chen in \cite{ChenChaos0} in the presence of external field (see also \cite{ChenChaos15} for a recent result that covers both cases). 

Although the statement of chaos in temperature in (\ref{Th1Eq}) looks very simple and does not refer to the properties of the two individual systems explicitly, the only known proof at the moment presented below passes through the entire Parisi ansatz and utilizes both ultrametricity/clustering of the overlaps and formula for the free energy. In fact, we will generalize the proof of ultrametricity in \cite{PUltra} and obtain \emph{joint clustering} for coupled systems at different temperatures in Theorem \ref{ThCluster} below. Joint clustering can be viewed as a kind of symmetry between the two systems, because it implies that a neighbourhood inside one system coincides with the neighbourhood of the same size of any nearby point from the second system (this property is expressed in the equations (\ref{Usymmetry})). It will be proved using a new invariance property for coupled systems analogous to the invariance property from the proof of ultrametricity for one system in \cite{PUltra}. As in \cite{PUltra}, the new invariance property is based on the strong form of the Ghirlanda-Guerra identities \cite{GuerraGG, GG}, which itself was obtained in \cite{PGGmixed} as a consequence of the Parisi formula for the free energy. In addition to yielding joint clustering, this invariance property possesses a certain built-in asymmetry due to the fact that two temperatures are different. This asymmetry turns out to be incompatible with the symmetry expressed by the joint clustering, under a certain assumption on the distributions of the overlaps within the two system which will be called \emph{uncoupled systems}, and will allow us to rule out `large' values of the cross-overlap. Another argument will rule out `intermediate' values.

In the complementary \emph{coupled} case, the two systems conspire to hide this asymmetry on some non-trivial interval of possible values of the cross-overlap, in some sense. Somewhat miraculously, this case turns out to be perfectly suited for another well-known approach based on Guerra's replica symmetry breaking interpolation \cite{Guerra}, proposed by Talagrand after his original proof of the Parisi formula in \cite{TPF}. Although very natural, so far this approach has not been used successfully on its own to prove ultrametricity or chaos in temperature in any general case because of seemingly intractable technical difficulties, and the case of the coupled systems that arises in this paper is, perhaps, the only non-trivial known example where these obstacles spontaneously disappear. We emphasize that this step also relies  on the validity of the Parisi formula for the free energy.

Combining these two very different types of techniques will allow us to eliminate all `large' and `intermediate' values of the cross-overlap. The proof will then be concluded by appealing to the results of Chen in \cite{ChenChaos} which handled `small' values of the overlaps in full generality.

The new invariance property will be proved and used in the infinite-volume limit $N\to\infty$. For this purpose, we will first need to construct an object that will play the role of asymptotic Gibbs measures for coupled systems, which will be done in the next section. The construction will be based on a version of the Dovbysh-Sudakov representation \cite{DS} for coupled systems.

\section{Asymptotic Gibbs measures for coupled systems}

This section may be skipped at first reading and Theorem \ref{ThDSpair} can be used as a black box, because the techniques in its proof are not directly related to the rest of the paper.

Consider a pair $G_N^1, G_N^2$ of random probability measures on $\{-1,+1\}^N$ that are not necessarily independent. The main example we have in mind are, of course, the Gibbs measures above. In this section, it will be convenient to denote replicas from both measures by $\sigma^\ell$ but over different sets of indices, $\ell\leq 0$ and $\ell\geq 1$, so we let $(\sigma^\ell)_{\ell\leq 0}$ be an i.i.d. sequence of replicas from $G_N^1$ and $(\sigma^\ell)_{\ell\geq 1}$ be an i.i.d. sequence of replicas from $G_N^2$. Let 
\begin{equation}
R^N = \bigl(R^N_{\ell,\ell'}\bigr)_{\ell,\ell'\in\mathbb{Z}}
= \bigl(\tilde{\sigma}^{\ell}\cdot \tilde{\sigma}^{\ell'}\bigr)_{\ell,\ell'\in\mathbb{Z}}
\label{arrayRN}
\end{equation} 
be the array of all their overlaps. Suppose that $R^N$ converges in distributions under $\e (G_N^1\times G_N^2)^{\otimes\infty}$ to some array $R = (R_{\ell,\ell'})_{\ell,\ell'\in\mathbb{Z}}$. Notice that this array is symmetric nonnegative-definite and also (partially) weakly exchangeable,
\begin{equation}
\bigl(R_{\pi(\ell),\pi(\ell')}\bigr)_{\ell,\ell'\in\mathbb{Z}}
\stackrel{d}{=}
\bigl(R_{\ell,\ell'}\bigr)_{\ell,\ell'\in\mathbb{Z}},
\label{permuteeq}
\end{equation}
but \emph{not} under all permutations of integers $\mathbb{Z}$, but only those permutations that map positive integers into positive and non-positive into non-positive. It turns out that, as a result, the off-diagonal elements are again generated by a pair of random measures which in the thermodynamic limit live on a separable Hilbert space. This is the analogue of the Dovbysh-Sudakov representation \cite{DS}.
\begin{theorem}\label{ThDSpair}
There exists a pair of random measures $G_1$ and $G_2$ on a separable Hilbert space $H$ (not necessarily independent) such that
$$
\bigl(R_{\ell,\ell'}\bigr)_{\ell\not =\ell'\in\mathbb{Z}}
\stackrel{d}{=}
\bigl(\sigma^{\ell}\cdot \sigma^{\ell'}\bigr)_{\ell\not =\ell'\in\mathbb{Z}},
$$
where $(\sigma^\ell)_{\ell\leq 0}$ is an i.i.d. sample from $G_1$ and $(\sigma^\ell)_{\ell\geq 1}$ is an i.i.d. sample from $G_2$. 
\end{theorem}
This will be a consequence of the following analogue of the Aldous-Hoover representation \cite{Aldous,Hoover}. Consider a pair of random arrays  
\begin{equation}
\bigl(s^1_{\ell,\ell'}\bigr)_{\ell,\ell'\geq 1}
\,\,\mbox{ and }\,\, 
\bigl(s^2_{\ell,\ell'}\bigr)_{\ell,\ell'\geq 1}
\label{arrayss}
\end{equation}
that are separately exchangeable in the first coordinate and jointly exchangeable in the second coordinate, that is,
\begin{equation}
\Bigl(\bigl(s^1_{\pi_1(\ell),\rho(\ell')}\bigr)_{\ell,\ell'\geq 1}, \bigl(s^2_{\pi_2(\ell),\rho(\ell')}\bigr)_{\ell,\ell'\geq 1}\Bigr)
\stackrel{d}{=}
\Bigl(\bigl(s^1_{\ell,\ell'}\bigr)_{\ell,\ell'\geq 1}, \bigl(s^2_{\ell,\ell'}\bigr)_{\ell,\ell'\geq 1}\Bigr)
\label{separatelyexch}
\end{equation}
for any permutations $\pi_1,\pi_2, \rho$ of finitely many coordinates. Then the following holds.
\begin{theorem}\label{ThAH2}
If (\ref{separatelyexch}) holds then there exist two measurable functions $\sigma_1,\sigma_2:[0,1]^4\to\Reals$ such that the arrays in (\ref{arrayss}) can be generated in distribution by
$$
\Bigl(\bigl(s^1_{\ell,\ell'}\bigr)_{\ell,\ell'\geq 1}, \bigl(s^2_{\ell,\ell'}\bigr)_{\ell,\ell'\geq 1}\Bigr)
\stackrel{d}{=}
\Bigl(\bigl(\sigma_1(w, u_\ell^1, v_{\ell'},x_{\ell,\ell'}^1)\bigr)_{\ell,\ell'\geq 1}, \bigl( \sigma_2(w, u_\ell^2, v_{\ell'},x_{\ell,\ell'}^2)\bigr)_{\ell,\ell'\geq 1}\Bigr),
$$
where all the arguments are i.i.d. uniform random variables on $[0,1]$.
\end{theorem}
The proofs of both theorems will be a simple modification of the arguments in Austin \cite{Austin}, and are based on the following version of de Finetti's theorem. Suppose that a random sequence $(s_\ell)_{\ell\geq 1}$ and random element $Z$ take values in some complete separable metric spaces $(S,{\cal S})$ and $(S',{\cal S}')$ with the Borel $\sigma$-algebras, and suppose that
\begin{eqnarray}
\bigl(Z,(s_\ell)_{\ell\geq 1}\bigr) \stackrel{d}{=} \bigl(Z,(s_{\pi(\ell)} \bigr)_{\ell\geq 1})
\label{eq:cond-exch}
\end{eqnarray}
for any permutation $\pi$ of finitely many coordinates. In this case, the almost sure limit
\begin{equation}
\eta = \lim_{n\to\infty} \frac{1}{n}\sum_{\ell=1}^n \delta_{s_\ell}
\label{empmeas}
\end{equation}
in the space of probability measures on $(S,{\cal S})$ (with the topology of weak convergence) is called the empirical measure of the sequence $(s_\ell)_{\ell\geq 1}$, if it exists. The following holds (see Proposition 1.4, Corollary 1.5 and Corollary 1.6 from \cite{Kallenberg}).
\begin{theorem}\label{ThdeFH}
If (\ref{eq:cond-exch}) holds then the empirical measure (\ref{empmeas}) exists almost surely and, conditionally on $\eta$, the sequence $(s_\ell)_{\ell\geq 1}$ is i.i.d. with the distribution $\eta$ and independent of $Z$.
\end{theorem}
The main part of the proof of Theorem \ref{ThAH2} will be based on the following observation. Suppose that we have two random sequences of pairs $(t^1_\ell,s^1_\ell)_{\ell\geq 1}$ and $(t^2_\ell,s^2_\ell)_{\ell\geq 1}$ with coordinates in complete separable metric spaces, which are separately exchangeable,
\begin{eqnarray}
\Bigl( 
\bigl(t^1_\ell, s^1_\ell\bigr)_{\ell\geq 1} , \bigl(t^2_\ell, s^2_\ell\bigr)_{\ell\geq 1} 
\Bigr)
\stackrel{d}{=} 
\Bigl(
\bigl(t^1_{\pi(\ell)}, s^1_{\pi(\ell)}\bigr)_{\ell\geq 1}, \bigl(t^2_{\rho(\ell)}, s^2_{\rho(\ell)}\bigr)_{\ell\geq 1}
\Bigr)
\end{eqnarray}
for any permutations $\pi$ and $\rho$ of finitely many coordinates. Consider the empirical measures
$$
\eta^1 =  \lim_{n\to\infty} \frac{1}{n}\sum_{\ell=1}^n \delta_{(t^1_\ell,s^1_\ell)},\,\,
\eta_1^1 =  \lim_{n\to\infty} \frac{1}{n}\sum_{\ell=1}^n \delta_{t^1_\ell}.
$$ 
and
$$
\eta^2 =  \lim_{n\to\infty} \frac{1}{n}\sum_{\ell=1}^n \delta_{(t^2_\ell,s^2_\ell)},\,\,
\eta_1^2 =  \lim_{n\to\infty} \frac{1}{n}\sum_{\ell=1}^n \delta_{t^2_\ell}.
$$ 
Obviously, $\eta_1^j$ is the marginal of $\eta^j$ on the first coordinate and, moreover, $\eta_1^j$ is a measurable function of the sequence $t^j = (t^j_\ell)_{\ell\geq 1}$. Let us denote 
\begin{equation}
\mbox{$\eta = (\eta^1,\eta^2)$, $\eta_1 = (\eta_1^1,\eta_1^2)$ and $t=(t^1,t^2)$}.
\label{triple}
\end{equation}
The following holds.
\begin{lemma}\label{LemEXclass2}
Conditionally on $(t^1_\ell)_{\ell\geq 1}$ and $(t^2_\ell)_{\ell\geq 1}$, we can generate sequences $(s^1_\ell)_{\ell\geq 1}$ and $(s^2_\ell)_{\ell\geq 1}$ in distribution as 
$$
\Bigl(\bigl(s^1_\ell\bigr)_{\ell\geq 1}, \bigl(s^2_\ell\bigr)_{\ell\geq 1}\Bigr) 
\stackrel{d}{=} 
\Bigl(
\bigl(f_1(\eta_1, t_\ell^1, v, x^1_\ell)\bigr)_{\ell\geq 1},
\bigl(f_2(\eta_1, t_\ell^2, v,x^2_\ell)\bigr)_{\ell\geq 1}
\Bigr)
$$
for some measurable functions $f_1$ and $f_2$ and i.i.d. uniform random variables $v$, $(x^1_\ell)_{\ell\geq 1}$ and $(x^2_\ell)_{\ell\geq 1}$ on $[0,1]$.
\end{lemma}
\textbf{Proof.}
First, let us note how to generate the empirical measures $(\eta^1,\eta^2)$ given the sequences $t^1=(t^1_\ell)_{\ell\geq 1}$ and $t^2=(t^2_\ell)_{\ell\geq 1}$.  Since $\eta_1^j$ is the first marginal of $\eta^j$, $(t^j_\ell)_{\ell\geq 1}$ are i.i.d. from $\eta^j_1$. This means that, if we consider the triple $(\eta,\eta_1,t)$ defined in (\ref{triple}) then the conditional distribution of $t$ given $(\eta,\eta_1)$ depends only on $\eta_1$,
$$
\p\bigl(t\in\cdot\, \,\bigr|\, \eta,\eta_1)=\p\bigl(t\in\cdot\, \,\bigr|\,\eta_1).
$$ 
This means that $t$ and $\eta$ are independent given $\eta_1$ and, therefore,
$$
\p\bigl(\eta \in\cdot\, \,\bigr|\, t,\eta_1)=\p\bigl(\eta \in\cdot\, \,\bigr|\,\eta_1).
$$ 
On the other hand, $\eta_1$ is a function of $t$, so
$
\p\bigl(\eta\in\cdot\, \,\bigr|\, t,\eta_1)=\p\bigl(\eta\in\cdot\, \,\bigr|\, t)
$
and, thus,
$$
\p\bigl(\eta\in\cdot\, \,\bigr|\, t)=
\p\bigl(\eta\in\cdot\, \,\bigr|\, \eta_1).
$$
In other words, to generate $\eta$ given $t$, we can simply compute $\eta_1$ and generate  $\eta$ given $\eta_1$. This means that, by the standard coding in terms of uniform random variables on $[0,1]$ (see e.g. Lemma 1.4 in \cite{SKmodel}), we can generate $\eta = g(\eta_1,v)$ in distribution as a function of $\eta_1$ and independent uniform random variable $v$ on $[0,1]$. 

Since the sequences $(t^1_\ell,s^1_\ell)_{\ell\geq 1}$ and $(t^2_\ell,s^2_\ell)_{\ell\geq 1}$ are separately exchangeable, by Theorem \ref{ThdeFH}, conditionally on $(\eta^1,\eta^2)$, the sequence $(t^1_\ell,s^1_\ell)_{\ell\geq 1}$ is i.i.d. from $\eta^1$, sequence $(t^2_\ell,s^2_\ell)_{\ell\geq 1}$ is i.i.d. from $\eta^2$, and these two sequences are independent of each other.  Therefore, given $t$ and $\eta$, we can simply generate $s_\ell^j$ from the conditional distribution $\eta^j(s^j_\ell \in\cdot\, \,|\, t^j_\ell)$ independently over $\ell$ and $j$. Again, this means that we can generate $s^j_\ell = h_j(\eta^j,t^j_\ell,x^j_\ell)$ as a function of i.i.d. uniform random variables $x^j_\ell$ on $[0,1]$. Finally, recalling that $\eta^j=g^j(\eta_1,v)$, we can write
$$
s_\ell^j = h_j(g_j(\eta_1,v),t^j_\ell,x^j_\ell) = f_j(\eta_1, t^j_\ell, v,x^j_\ell)
$$
for some functions $f_j$. This finishes the proof.
\qed

\medskip
\medskip
\noindent
\textbf{Proof of Theorem \ref{ThAH2}.} Let us for convenience index the arrays $s^j_{\ell,\ell'}$ by $\ell\geq 1$ and $\ell'\in\mathbb{Z}$ instead of $\ell'\geq 1$. Let us denote 
$$
s_{\ell,\ell'}=\bigl(s^1_{\ell,\ell'}, s^2_{\ell,\ell'} \bigr),\,\,
X_{\ell'} = \bigl(s_{\ell,\ell'}\bigr)_{\ell\geq 1}
\,\,\mbox{ and }\,\,
X=\bigl(X_{\ell'}\bigr)_{\ell'\leq 0}.
$$ 
Since the sequence of columns $(X_{\ell'})_{\ell'\in\mathbb{Z}}$ is exchangeable and the empirical measure is a function of $X$, conditionally on $X$, the columns $(X_{\ell'})_{\ell'\geq 1}$ in the `right half' of the array are i.i.d.. If we describe the distribution of one column $X_1$ given $X$ then we can generate all columns $(X_{\ell'})_{\ell'\geq 1}$ independently from this distribution. Hence, our strategy will be to describe the distribution of $X_1$ given $X$, and then combine it with the structure of the distribution of $X$. Both steps will use exchangeability with respect to permutations of rows, because so far we have only used exchangeability with respect to permutations of columns. Let us denote 
$$
Y_\ell^j = \bigl(s^j_{\ell,\ell'} \bigr)_{\ell'\leq 0}.
$$
We want to describe the distribution of $X_1 = (s^1_{\ell,1}, s_{\ell,1}^2)_{\ell\geq 1}$ given $X = (Y_\ell^1,Y_\ell^2)_{\ell\geq 1}$ and we will use the fact that the sequences $(Y^1_\ell,s^1_{\ell,1})_{\ell\geq 1}$ and $(Y^2_\ell,s^2_{\ell,1})_{\ell\geq 1}$ are separately exchangeable. By Lemma \ref{LemEXclass2}, conditionally on $X=(Y^1_\ell, Y_\ell^2)_{\ell\geq 1}$, $(s^1_{\ell,1})_{\ell\geq 1}$ and $(s^1_{\ell,1})_{\ell\geq 1}$ can be generated in distribution as
$$
s^j_{\ell,1} = f_j(\eta_1, Y^j_\ell, v_1,x^j_{\ell,1}),
$$
where $\eta_1 = (\eta_1^1,\eta_1^2)$ and $\eta_1^j$ is the empirical measure of $(Y^j_\ell)$, and where instead of $v$ and $(x^j_\ell)$ we wrote $v_1$ and $(x^j_{\ell,1})$ to emphasize the first column index $1$. Since, conditionally on $X$, the columns $(X_{\ell'})_{\ell'\geq 1}$ in the `right half' of the array are i.i.d., we can generate
$$
s^j_{\ell,\ell'} = f_j(\eta_1, Y^j_\ell, v_{\ell'},x^j_{\ell,\ell'}),
$$
where $v_{\ell'}$ and $x^j_{\ell,\ell'}$ are i.i.d. uniform random variables on $[0,1]$. Finally, we use that $(Y_\ell^1)_{\ell\geq 1}$ and $(Y_\ell^2)_{\ell\geq 1}$ are separately exchangeable. By Lemma \ref{LemEXclass2}, conditionally on $\eta_1 = (\eta_1^1,\eta_1^2)$, $(Y^j_\ell)_{\ell\geq 1}$ are i.i.d. from $\eta_1^j$ and independent of each other. Again, coding in terms of uniform random variables on $[0,1]$, we can generate $\eta_1= h(w)$ as a function of a uniform random variable $w$ on $[0,1]$ and then generate $Y^j_\ell = Y_j(\eta_1,u^j_\ell)=Y_j(h(w),u^j_\ell)$ as functions of $\eta_1$ and i.i.d. uniform random variables $u^j_\ell$ on $[0,1]$. Plugging these into $f_j$ above gives 
$$
s^j_{\ell,\ell'} = \sigma_j(w, u^j_\ell, v_{\ell'},x^j_{\ell,\ell'})
$$ 
for some functions $\sigma_1$ and $\sigma_2$, which finishes the proof.
\qed

\medskip
\medskip
\noindent
\textbf{Proof of Theorem \ref{ThDSpair}.}
Since the array $R$ is nonnegative-definite, conditionally on $R$, we can generate a Gaussian vector $g$ in $\Reals^{\mathbb{Z}}$ with the covariance equal to $R$. Now, also conditionally on $R$, let $(g_i)_{i\geq 1}$ be independent copies of $g$. For each $i\geq 1$, let us denote the coordinates of $g_i$ by $g_{\ell,i}$ for $\ell \in \mathbb{Z}$. Then, since the array $R=(R_{\ell,\ell'})_{\ell,\ell'\geq 1}$ satisfies (\ref{permuteeq}) under all permutations of integers $\mathbb{Z}$ that map positive integers into positive and non-positive into non-positive, we have
\begin{equation}
\Bigl(\bigl(g_{\pi_1(\ell),\rho(i)}\bigr)_{\ell\geq 1,i\geq 1}, \bigl(g_{\pi_2(\ell),\rho(i)}\bigr)_{\ell\leq 0,i\geq 1}\Bigr)
\stackrel{d}{=}
\Bigl(\bigl(g_{\ell,i}\bigr)_{\ell\geq 1,i\geq 1}, \bigl(g_{\ell,i}\bigr)_{\ell\leq 0, i\geq 1}\Bigr)
\end{equation}
for any permutations $\pi_1,\pi_2, \rho$ of finitely many coordinates. This is precisely the property in (\ref{separatelyexch}), only here instead of using superscripts $1$ and $2$ we used different sets of subscripts $\ell\geq 1$ and $\ell\leq 0$. By Theorem \ref{ThAH2}, there exist two measurable functions $\sigma_1,\sigma_2:[0,1]^4\to\Reals$ such that these arrays can be generated in distribution by
$$
\Bigl(\bigl(g_{\ell,i}\bigr)_{\ell\geq 1,i\geq 1}, \bigl(g_{\ell,i}\bigr)_{\ell\leq 0, i\geq 1}\Bigr)
\stackrel{d}{=}
\Bigl(\bigl(\sigma_1(w, u_\ell, v_{i},x_{\ell,i}) \bigr)_{\ell\geq 1,i\geq 1},
\bigl(\sigma_2(w, u_\ell, v_{i},x_{\ell,i}) \bigr)_{\ell\leq 0, i\geq 1}\Bigr).
$$
By the strong law of large numbers (applied conditionally on $R$), for any $\ell\not  = \ell'$,
$$
\frac{1}{n} \sum_{i=1}^n g_{\ell,i} g_{\ell',i} \to R_{\ell,\ell'}
$$ 
almost surely as $n\to\infty$. Similarly, by the strong law of large numbers (now applied conditionally on $w$, $(u_{\ell})_{\ell\geq 1}$ and $(u_{\ell})_{\ell\leq 0}$), for any $\ell\not  = \ell'$,
$$
\frac{1}{n} \sum_{i=1}^n \sigma_j(w,u_{\ell},v_{i},x_{\ell,i}) \sigma_{j'}(w,u_{\ell'},v_{i},x_{\ell',i}) \to 
\e' \sigma_j(w,u_{\ell},v,x_1) \sigma_{j'}(w,u_{\ell'},v,x_{2})
$$ 
almost surely, where $\e'$ denotes the expectation with respect to the random variables $v, x_1, x_2$. Here $j=1$ if $\ell\geq 1$, $j=2$ if $\ell\leq 0$ and, similarly, for $j'$ and $\ell'$. Therefore, we showed  that
\begin{equation}
\bigl( R_{\ell,\ell'}\bigr)_{\ell \not = \ell'}
\stackrel{d}{=}
\bigl(\e' \sigma_j(w,u_{\ell},v,x_{1}) \sigma_{j'}(w,u_{\ell'},v,x_{2})\bigr)_{\ell\not =\ell'},
\label{chApp23DSdone}
\end{equation}
where $j=1$ if $\ell\geq 1$, $j=2$ if $\ell\leq 0$ and, similarly, for $j'$ and $\ell'$. If we denote
$$
\xoverline{\sigma}_{1}(w,u,v) = \int\! \sigma_1(w,u,v,x) \, dx,\,\,\,
\xoverline{\sigma}_{2}(w,u,v) = \int\! \sigma_2(w,u,v,x) \, dx
$$
then (\ref{chApp23DSdone}) can be rewritten as
\begin{equation}
\bigl( R_{\ell,\ell'}\bigr)_{\ell \not = \ell'}
\stackrel{d}{=}
\bigl(\e' \xoverline{\sigma}_j(w,u_{\ell},v) \xoverline{\sigma}_{j'}(w,u_{\ell'},v)\bigr)_{\ell\not =\ell'}.
\label{chApp23DSdone2}
\end{equation}
Notice that, for almost all $w$ and $u$, the functions $v\to \xoverline{\sigma}_{j}(w,u,v)$ are in $H = L^2([0,1], dv)$. Therefore, if we denote $\sigma^{\ell} = \xoverline{\sigma}_{j}(w,u_{\ell},\,\cdot\,),$ where $j=1$ if $\ell\geq 1$, $j=2$ if $\ell\leq 0$, then (\ref{chApp23DSdone2}) becomes
\begin{equation}
\bigl( R_{\ell,\ell'}\bigr)_{\ell\not = \ell'}
\stackrel{d}{=}
\bigl( \sigma^{\ell}\cdot \sigma^{\ell'} \bigr)_{\ell\not = \ell'}.
\label{chApp23DSdone3}
\end{equation}
It remains to observe that $(\sigma^{\ell})_{\ell\geq 1}$ is an i.i.d. sequence from the random measure $G_1$ on $H$ given by the image of the Lebesgue measure $du$ on $[0,1]$ by the map $u \to \xoverline{\sigma}_{1}(w,u,\,\cdot\,)$ and $(\sigma^{\ell})_{\ell\leq 0}$ is an i.i.d. sequence from the random measure $G_2$ on $H$ given by the image of the Lebesgue measure $du$ on $[0,1]$ by the map $u \to \xoverline{\sigma}_{2}(w,u,\,\cdot\,)$. This finishes the proof.
\qed

\section{Invariance properties}\label{SecInvariance}

In this section, we will show that if one starts with the Gibbs measures $G_{N}^1$ and $G_{N}^2$, considers an array $R = (R_{\ell,\ell'})_{\ell,\ell'\in\mathbb{Z}}$ given by any subsequential limit in distribution of the overlap arrays $R^N$ in (\ref{arrayRN}), and then considers a pair of random measures $G_1$ and $G_2$ constructed in Theorem \ref{ThDSpair}, then these measures will satisfy certain joint invariance properties, which will be consequences of the Ghirlanda-Guerra identities \cite{GuerraGG, GG}.

The starting point is the strong form of the Ghirlanda-Guerra identities  proved in \cite{PGGmixed}. In the form of the concentration of Hamiltonian, these identities say that
\begin{equation}
\lim_{N\to\infty} \frac{1}{N}\e 
\bigl\la
\bigl|
H_{N,p}(\sigma) - \e\la H_{N,p}(\sigma)\ra
\bigr|
\bigr\ra
=0
\label{GGHamp}
\end{equation}
whenever $\gamma_p\not = 0$ in (\ref{Hammixed}) and, similarly,
\begin{equation}
\lim_{N\to\infty} \frac{1}{N}\e 
\bigl\la
\bigl|
H_{N,p}(\rho) - \e\la H_{N,p}(\rho)\ra
\bigr|
\bigr\ra
=0.
\label{GGHamp2}
\end{equation}
In the rest of the paper we will denote the ratio of inverse temperature parameters by
\begin{equation}
\kappa = \frac{\beta_1}{\beta_2} \not=1.
\end{equation}
Let us consider a continuous bounded function $\Phi$ of the overlaps
$$
\tilde{\sigma}^\ell\cdot \tilde{\sigma}^{\ell'},
\tilde{\sigma}^\ell\cdot \tilde{\rho}^{\ell'},
\tilde{\rho}^\ell\cdot \tilde{\rho}^{\ell'}
\mbox{ for $\ell,\ell'\leq n$}
$$
of $n$ replicas $(\sigma^\ell,\rho^\ell)_{\ell\leq n}$ from $G_{N}^1\times G_{N}^2$. If we use (\ref{GGHamp}) to write
$$
\frac{1}{N}
\e 
\bigl\la
\Phi H_{N,p}(\sigma^1)
\bigr\ra
\approx
\frac{1}{N}
\e\bigl \la H_{N,p}(\sigma^1)\bigr\ra
\e\bigl \la \Phi\bigr\ra
$$
and then use Gaussian integration by parts on both sides (see e.g. Lemma 1.1 in \cite{SKmodel}), we get
\begin{align*}
\e\bigl \la \Phi\bigl(\beta_1(\tilde{\sigma}^1\cdot \tilde{\sigma}^{n+1})^p+\beta_2(\tilde{\sigma}^1\cdot \tilde{\rho}^{n+1})^p \bigr)\bigr\ra
& \approx
\frac{1}{n} \sum_{\ell=1}^n \beta_2 \e\bigl \la \Phi (\tilde{\sigma}^1\cdot \tilde{\rho}^{\ell})^p \bigr\ra
\\
& + \frac{1}{n}\beta_1 \e\bigl \la \Phi\bigr\ra \e\bigl \la (\tilde{\sigma}^1\cdot \tilde{\sigma}^{2})^p \bigr\ra +
\frac{1}{n} \sum_{\ell=2}^n \beta_1 \e\bigl \la \Phi (\tilde{\sigma}^1\cdot \tilde{\sigma}^{\ell})^p \bigr\ra
\end{align*}
or, equivalently,
\begin{align}
\e\bigl \la \Phi\bigl((\tilde{\sigma}^1\cdot \tilde{\sigma}^{n+1})^p+\frac{1}{\kappa}(\tilde{\sigma}^1\cdot \tilde{\rho}^{n+1})^p \bigr)\bigr\ra
& \approx 
\frac{1}{n} \sum_{\ell=1}^n \frac{1}{\kappa} \e\bigl \la \Phi (\tilde{\sigma}^1\cdot \tilde{\rho}^{\ell})^p \bigr\ra
\label{GGN1}
\\
&+  \frac{1}{n} \e\bigl \la \Phi\bigr\ra \e\bigl \la (\tilde{\sigma}^1\cdot \tilde{\sigma}^{2})^p \bigr\ra
+
\frac{1}{n} \sum_{\ell=2}^n  \e\bigl \la \Phi (\tilde{\sigma}^1\cdot \tilde{\sigma}^{\ell})^p \bigr\ra.
\nonumber
\end{align}
Similarly, using  (\ref{GGHamp}), we can write
\begin{align}
\e\bigl \la \Phi\bigl((\tilde{\rho}^1\cdot \tilde{\rho}^{n+1})^p+\kappa(\tilde{\rho}^1\cdot \tilde{\sigma}^{n+1})^p \bigr)\bigr\ra
\approx &
\frac{1}{n} \sum_{\ell=1}^n \kappa \e\bigl \la \Phi (\tilde{\rho}^1\cdot \tilde{\sigma}^{\ell})^p \bigr\ra
\label{GGN2}
\\
& + \frac{1}{n} \e\bigl \la \Phi\bigr\ra \e\bigl \la (\tilde{\rho}^1\cdot \tilde{\rho}^{2})^p \bigr\ra
+
\frac{1}{n} \sum_{\ell=2}^n  \e\bigl \la \Phi (\tilde{\rho}^1\cdot \tilde{\rho}^{\ell})^p \bigr\ra.
\nonumber
\end{align}
If we consider any limit of the array of all overlaps in distribution, in the limit, by Theorem \ref{ThDSpair}, the overlaps can be generated from some random measure $G_1\times G_2$ on $H^2$, where $H$ is a separable Hilbert space. If we denote by $(\sigma^\ell,\rho^{\ell})_{\ell\geq 1}$ i.i.d. replicas from $G_1\times G_2$ and if we continue to use the notation $\la\,\cdot\,\ra$ for the average with respect to $(G_1\times G_2)^{\otimes \infty}$ then the above approximate identities will become exact identities
\begin{align*}
\e\bigl \la \Phi\bigl(({\sigma}^1\cdot {\sigma}^{n+1})^p+\frac{1}{\kappa}({\sigma}^1\cdot {\rho}^{n+1})^p \bigr)\bigr\ra
= &
\frac{1}{n} \sum_{\ell=1}^n \frac{1}{\kappa} \e\bigl \la \Phi ({\sigma}^1\cdot {\rho}^{\ell})^p \bigr\ra
\\
& + \frac{1}{n} \e\bigl \la \Phi\bigr\ra \e\bigl \la ({\sigma}^1\cdot {\sigma}^{2})^p \bigr\ra
+
\frac{1}{n} \sum_{\ell=2}^n  \e\bigl \la \Phi ({\sigma}^1\cdot {\sigma}^{\ell})^p \bigr\ra
\end{align*}
and
\begin{align*}
\e\bigl \la \Phi\bigl(({\rho}^1\cdot {\rho}^{n+1})^p+\kappa({\rho}^1\cdot {\sigma}^{n+1})^p \bigr)\bigr\ra
= &
\frac{1}{n} \sum_{\ell=1}^n \kappa \e\bigl \la \Phi ({\rho}^1\cdot {\sigma}^{\ell})^p \bigr\ra
\\
& + \frac{1}{n} \e\bigl \la \Phi\bigr\ra \e\bigl \la ({\rho}^1\cdot {\rho}^{2})^p \bigr\ra
+
\frac{1}{n} \sum_{\ell=2}^n  \e\bigl \la \Phi ({\rho}^1\cdot {\rho}^{\ell})^p \bigr\ra,
\end{align*}
where $\Phi=\Phi(R^n)$ is now a continuous function of the Gram matrix $R^n$ of the overlaps
$$
{\sigma}^\ell\cdot {\sigma}^{\ell'},
{\sigma}^\ell\cdot {\rho}^{\ell'},
{\rho}^\ell\cdot {\rho}^{\ell'}
\mbox{ for $\ell,\ell'\leq n$}
$$
of $n$ replicas $(\sigma^\ell,\rho^\ell)_{\ell\leq n}$. Since we work with generic models and linear span of functions $x^p$ for even $p\geq 2$ with $\gamma_p \not = 0$ is dense in $C([0,1],\|\,\cdot\,\|_\infty)$, we get that\begin{align}
\e\bigl \la \Phi\bigl(\psi({\sigma}^1\cdot {\sigma}^{n+1})+\frac{1}{\kappa}\psi({\sigma}^1\cdot {\rho}^{n+1}) \bigr)\bigr\ra
= &
\frac{1}{n} \sum_{\ell=1}^n \frac{1}{\kappa} \e\bigl \la \Phi \psi({\sigma}^1\cdot {\rho}^{\ell}) \bigr\ra
\label{GG1}
\\
& + \frac{1}{n} \e\bigl \la \Phi\bigr\ra \e\bigl \la \psi({\sigma}^1\cdot {\sigma}^{2}) \bigr\ra
+
\frac{1}{n} \sum_{\ell=2}^n  \e\bigl \la \Phi \psi({\sigma}^1\cdot {\sigma}^{\ell}) \bigr\ra
\nonumber
\end{align}
and
\begin{align}
\e\bigl \la \Phi\bigl(\psi({\rho}^1\cdot {\rho}^{n+1})+\kappa\psi({\rho}^1\cdot {\sigma}^{n+1}) \bigr)\bigr\ra
= &
\frac{1}{n} \sum_{\ell=1}^n \kappa \e\bigl \la \Phi \psi({\rho}^1\cdot {\sigma}^{\ell}) \bigr\ra
\label{GG2}
\\
& + \frac{1}{n} \e\bigl \la \Phi\bigr\ra \e\bigl \la \psi({\rho}^1\cdot {\rho}^{2}) \bigr\ra
+
\frac{1}{n} \sum_{\ell=2}^n  \e\bigl \la \Phi \psi({\rho}^1\cdot {\rho}^{\ell}) \bigr\ra
\nonumber
\end{align}
for any continuous function $\psi$ on $[-1,1]$, which is symmetric,
\begin{equation}
\psi(x) = \psi(|x|) \mbox{ for } x\in[-1,1]. 
\end{equation}
Of course, the identities (\ref{GG1}) and (\ref{GG2}) then hold for measurable bounded functions $\Phi$ and $\psi$, if $\psi$ is symmetric. Such identities first appeared in Chen, Panchenko \cite{ChenChaos2} and later used by Chen in \cite{ChenChaos}, but here they will be used very differently, by way of the following invariance property.

Let us denote a generic point in $H$ by $\tau$, which could be either $\sigma$ or $\rho$. Given $n\geq 1$, we consider $n$ bounded measurable functions $f_1,\ldots, f_n$ on $[-1,1]$ that are symmetric, $f_j(x) = f_j(|x|),$ and let
$$
F(\tau,\sigma^1,\ldots,\sigma^n) = f_1(\tau\cdot\sigma^1)+\ldots+f_n(\tau\cdot\sigma^n).
$$
For $1\leq \ell\leq n$, we define
$$
F_{\ell}(\tau,\sigma^1,\ldots,\sigma^n) = F(\tau,\sigma^1,\ldots,\sigma^n)
 - f_{\ell}(\tau\cdot\sigma^{\ell})+ \e \bigl \la f_{\ell}\bigl(\sigma^1\cdot \sigma^2 \bigr) \bigr\ra.
$$
Similarly, let us consider $n$ bounded measurable functions $g_1,\ldots, g_n$ symmetric on $[-1,1]$ and let
$$
G(\tau,\rho^1,\ldots,\rho^n) = g_1(\tau\cdot\rho^1)+\ldots+g_n(\tau\cdot\rho^n)
$$
and, for $1\leq \ell\leq n$,
$$
G_{\ell}(\tau,\rho^1,\ldots,\rho^n) = G(\tau,\rho^1,\ldots,\rho^n) 
 - g_{\ell}(\tau \cdot\rho^{\ell})+ \e \bigl \la g_{\ell}\bigl(\rho^1\cdot\rho^2 \bigr) \bigr\ra.
$$
For $\ell\leq n$, we define
\begin{align*}
D_\ell(\sigma,\rho)
& = F_\ell(\sigma,\sigma^1,\ldots,\sigma^n) + \kappa G(\sigma,\rho^1,\ldots,\rho^n)
\\
& + \frac{1}{\kappa} F(\rho,\sigma^1,\ldots,\sigma^n) + G_\ell(\rho,\rho^1,\ldots,\rho^n)
\end{align*}
and, for $\ell\geq n$, we define $D_\ell(\sigma,\rho)=D(\sigma,\rho)$, where
\begin{align*}
D(\sigma,\rho)
& = F(\sigma,\sigma^1,\ldots,\sigma^n) + \kappa G(\sigma,\rho^1,\ldots,\rho^n)
\\
& + \frac{1}{\kappa} F(\rho,\sigma^1,\ldots,\sigma^n) + G(\rho,\rho^1,\ldots,\rho^n).
\end{align*}
In both cases, for simplicity of notation, we omit the dependence on $(\sigma^\ell,\rho^\ell)_{\ell\leq n}$. The following holds.
\begin{theorem}\label{ch45Th1}
Let $\Phi$ be a bounded measurable function of the overlaps $R^n$ of $n$ replicas. Then,
\begin{equation}
\e \bigl\la  \Phi(R^n) \bigr\ra
=
\e\biggl\la
\frac{ \Phi(R^n) \exp \sum_{\ell=1}^n D_\ell(\sigma^\ell,\rho^\ell)}
{\bigl\la\exp D(\sigma,\rho)\bigr\ra_{\hspace{-0.3mm}\mathunderscore}^n}
\biggr\ra.
\label{ch45main}
\end{equation}
where the average $\la\,\cdot\,\ra_{\hspace{-0.3mm}\mathunderscore}$ with respect to $G_1\times G_2$ in the denominator is in $(\sigma,\rho)$ only for fixed $(\sigma^\ell,\rho^\ell)_{\ell\leq n}$, and the outside average $\la\,\cdot\, \ra$ of the ratio is in $(\sigma^\ell,\rho^\ell)_{\ell\leq n}$.
\end{theorem}
\noindent\textbf{Proof.}
Without loss of generality, let us assume that $|\Phi| \leq 1$ and suppose that $|f_{\ell}|\leq L$ and $|g_{\ell}|\leq L$ for all $\ell\leq n$ for some large enough $L.$ For $t\geq 0$, let us define
\begin{equation}
\varphi(t) = 
\e\biggl\la
\frac{ \Phi \exp \sum_{\ell=1}^n tD_\ell(\sigma^\ell,\rho^\ell)}
{\bigl\la\exp tD(\sigma,\rho)\bigr\ra_{\hspace{-0.3mm}\mathunderscore}^n}
\biggr\ra.
\label{ch45varphitdefine}
\end{equation}
We will show that the Ghirlanda-Guerra identities (\ref{GG1}) and (\ref{GG2}) imply that the function $\varphi(t)$ is constant for all $t\geq 0$, proving the statement of the theorem, $\varphi(0)=\varphi(1).$ For $k\geq 1$, let us denote
$$
\Pi_{n+k} = \sum_{\ell=1}^{n+k-1} D_\ell(\sigma^\ell,\rho^\ell) -(n+k-1)D_{n+k}(\sigma^{n+k},\rho^{n+k}).
$$
Using that the average $\la\,\cdot\,\ra_{\hspace{-0.3mm}\mathunderscore}$ is in $(\sigma,\rho)$ only, one can check by induction on $k$ that
$$
\varphi^{(k)}(t) = 
\e\biggl\la
\frac{ \Phi \Pi_{n+1}\cdots \Pi_{n+k} \exp \sum_{\ell=1}^{n+k} tD_\ell(\sigma^\ell,\rho^\ell)}
{\bigl\la\exp tD(\sigma,\rho)\bigr\ra_{\hspace{-0.3mm}\mathunderscore}^{n+k}}
\biggr\ra.
$$
Next, we will show that $\varphi^{(k)}(0)=0.$ If we introduce the notation 
$$
\Phi' = \Phi \Pi_{n+1}\cdots \Pi_{n+k-1},
$$ 
then $\Phi'$ is a function of the overlaps of $n+k-1$ replicas $(\sigma^\ell,\rho^\ell)_{\ell\leq n+k-1}$ and $\varphi^{(k)}(0)=\e\la \Phi' \Pi_{n+k}\ra.$ On the other hand, $\Pi_{n+k}$ can be written as
\begin{align*}
\sum_{j=1}^n \Bigl(
\sum_{\ell=1}^{n+k-1} \frac{1}{\kappa} f_j(\rho^\ell\cdot\sigma^j)
&+ \sum_{\ell\not = j,\ell=1}^{n+k-1}f_j(\sigma^\ell\cdot\sigma^j)
+ \e \bigl\la f_j(\sigma^1\cdot\sigma^2) \bigr\ra
\\
& -(n+k-1) f_{j}(\sigma^{n+k}\cdot\sigma^j)- (n+k-1) \frac{1}{\kappa} f_{j}(\rho^{n+k}\cdot\sigma^j)
\Bigr)
\\
+\sum_{j=1}^n \Bigl(
\sum_{\ell=1}^{n+k-1} \kappa g_j(\sigma^\ell\cdot\rho^j)
&+ \sum_{\ell\not = j,\ell=1}^{n+k-1}g_j(\rho^\ell\cdot\rho^j)
+ \e \bigl\la g_j(\rho^1\cdot\rho^2) \bigr\ra
\\
& -(n+k-1) g_{j}(\rho^{n+k}\cdot\rho^j)- (n+k-1) \kappa g_{j}(\sigma^{n+k}\cdot\rho^j)
\Bigr).
\end{align*}
Therefore, applying the Ghirlanda-Guerra identities (\ref{GG1}) or (\ref{GG2}) to each term $j\leq n$, we get that  
$$
\varphi^{(k)}(0)=\e\la \Phi' \Pi_{n+k}\ra = 0.
$$
If we denote $M=L(2+\kappa+\kappa^{-1})$ then $|D_{\ell}| \leq M n$ and $|\Pi_{n+k}| \leq 2Mn$. We can then bound
\begin{align*}
\bigl|\varphi^{(k)}(t) \bigr|
&\,\leq \,
\Bigl(\prod_{\ell=1}^{k} 2M (n+\ell-1)n \Bigr) \,
\e\biggl\la
\frac{ \Phi \exp \sum_{\ell=1}^{n+k} tD_\ell(\sigma^\ell,\rho^\ell)}
{\bigl\la\exp tD(\sigma,\rho)\bigr\ra_{\hspace{-0.3mm}\mathunderscore}^{n+k}}
\biggr\ra
\\
&\,=\,
\Bigl(\prod_{\ell=1}^{k} 2M (n+\ell-1)n\Bigr) \,
\e\biggl\la
\frac{ \Phi \exp \sum_{\ell=1}^{n} tD_\ell(\sigma^\ell,\rho^\ell)}
{\bigl\la\exp tD(\sigma,\rho)\bigr\ra_{\hspace{-0.3mm}\mathunderscore}^{n}}
\biggr\ra,
\end{align*}
where the equality in the second line follows from the fact that the denominator and $\Phi$ depend only on $(\sigma^\ell,\rho^\ell)_{\ell\leq n}$ and the average of the numerator in $(\sigma^{\ell},\rho^\ell)$ for each $n<\ell\leq n+k$ will cancel exactly one factor in the denominator. Moreover, if we consider an arbitrary $T>0$, using that $|D_{\ell}| \leq M n$, the last ratio can be bounded by $\exp (2M T n^2)$ for $0\leq t\leq T$ and, therefore,
$$
\max_{0\leq t\leq T}\bigl|\varphi^{(k)}(t) \bigr| \leq \exp( 2M T n^2) \frac{(n+k-1)!}{(n-1)!}\,  (2Mn)^k.
$$
Since we proved above that $\varphi^{(k)}(0) = 0$ for all $k\geq 1$, using Taylor's expansion, we can write
$$
\bigl|\varphi(t)-\varphi(0) \bigr| 
\leq 
 \max_{0\leq s\leq t} \frac{|\varphi^{(k)}(s)|}{k!}t^k
\leq
\exp(2M T n^2)  \frac{(n+k-1)! }{k! \,(n-1)!} (2Mn t)^k.
$$
Letting $k\to \infty$ proves that $\varphi(t)=\varphi(0)$ for $0\leq t<(2Mn)^{-1}.$ This implies that for any $t_0<(2Mn)^{-1}$ we have $\varphi^{(k)}(t_0)=0$ for all $k\geq 1$ and, again, by Taylor's expansion for $t_0\leq t \leq T,$
\begin{align*}
\bigl|\varphi(t)-\varphi(t_0)\bigr| 
&\leq 
\max_{t_0\leq s\leq t} \frac{|\varphi^{(k)}(s)|}{k!}(t-t_0)^k
\\
&\leq 
\exp(2M T n^2)  \frac{(n+k-1)! }{k! \,(n-1)!} \bigl(2Mn (t-t_0)\bigr)^k.
\end{align*}
Letting $k\to\infty$ proves that $\varphi(t) = \varphi(0)$ for $0\leq t< 2(2Mn)^{-1}.$ We can proceed in the same fashion to prove this equality for all $t< T$ and note that $T$ was arbitrary.
\qed

Next, consider a finite index set $\A$ and let $(B_\alpha)_{\alpha\in\A}$ be some partition of $H\times H$ such that, for each $\alpha\in\A$, the indicator $I_{B_\alpha}(\sigma,\rho) = I((\sigma,\rho)\in B_\alpha)$ 
is a measurable function of the overlaps $R^n$ and 
$$
\sigma\cdot\sigma^\ell, \sigma\cdot\rho^\ell, \rho\cdot\sigma^\ell, \rho\cdot\rho^\ell
\mbox{ for } \ell\leq n.
$$ 
In other words, the sets in the partition are expressed in terms of some conditions on the scalar products between $\sigma,\sigma^1,\ldots,\sigma^n,\rho,\rho^1,\ldots,\rho^n$. Let
\begin{equation}
W_\alpha=W_\alpha\bigl((\sigma^\ell,\rho^\ell)_{\ell\leq n}\bigr)=G(B_\alpha)
\label{ch45WA}
\end{equation}
be the weights of the sets in this partition with respect to the measure $G=G_1\times G_2$.
Let us define a map $T$ by
\begin{equation}
W=(W_\alpha)_{\alpha\in\A}\to T(W) = 
\biggl(\frac{\la I_{B_\alpha}(\sigma,\rho) \exp D(\sigma,\rho) \ra_{\hspace{-0.3mm}\mathunderscore}}
{\la \exp D(\sigma,\rho)\ra_{\hspace{-0.3mm}\mathunderscore} } \biggr)_{\alpha\in\A}.
\label{ch45TA}
\end{equation}
Then the following holds.
\begin{theorem}\label{ch45Th2}
Let $\varphi$ be a bounded measurable function of the overlaps $R^n$ of $n$ replicas and the weights $W$ in (\ref{ch45WA}). Then,
\begin{equation}
\e \bigl\la  \varphi(R^n,W) \bigr\ra
=
\e\biggl\la
\frac{ \varphi(R^n, T(W)) \exp \sum_{\ell=1}^n D_\ell(\sigma^\ell,\rho^\ell)}
{\bigl\la\exp D(\sigma,\rho)\bigr\ra_{\hspace{-0.3mm}\mathunderscore}^n}
\biggr\ra.
\label{ch45main2}
\end{equation}
where the average $\la\,\cdot\,\ra_{\hspace{-0.3mm}\mathunderscore}$ with respect to $G_1\times G_2$ in the denominator is in $(\sigma,\rho)$ only for fixed $(\sigma^\ell,\rho^\ell)_{\ell\leq n}$, and the outside average $\la\,\cdot\, \ra$ of the ratio is in $(\sigma^\ell,\rho^\ell)_{\ell\leq n}$.
\end{theorem}
\noindent\textbf{Proof.}
Let $n_\alpha\geq 0$ be some integers for $\alpha\in\A$ and let $m=n+\sum_{\alpha\in\A} n_\alpha.$ Let $(S_\alpha)_{\alpha\in\A}$ be any partition of $\{n+1,\ldots,m\}$ such that the cardinalities $|S_\alpha |=n_\alpha.$ Consider a continuous function $\Phi = \Phi(R^n)$ of the overlaps of $n$ replicas and let 
$$
\Phi' = \Phi(R^n) \prod_{\alpha\in\A}\varphi_\alpha,
\,\mbox{ where }\,
\varphi_\alpha = I \bigl((\sigma^{\ell},\rho^\ell) \in B_\alpha, \forall \ell\in S_\alpha \bigr).
$$
We will apply Theorem \ref{ch45Th1} to the function $\Phi'$, but since it now depends on $m$ coordinates,  we have to choose $2m$ bounded measurable functions $f_1,\ldots, f_m$ and  $g_1,\ldots, g_m$ in the definition of $D$ and $D_\ell$ above. We will choose the functions $f_1,\ldots, f_n$ and  $g_1,\ldots, g_n$ to be arbitrary and we let 
\begin{equation}
f_{n+1}=\ldots=f_m=g_{n+1}=\ldots = g_m = 0. 
\label{fgzeros}
\end{equation}
First of all, integrating out the coordinates $(\sigma^{\ell},\rho^\ell)_{\ell>n}$, the left hand side of (\ref{ch45main}) can be written  as
\begin{equation}
\e \bigl\la \Phi'  \bigr\ra 
=
\e \Bigl\la \Phi(R^n) \prod_{\alpha\in\A} \varphi_\alpha \Bigr\ra
=
\e \Bigl\la\Phi(R^n) \prod_{\alpha\in\A} W_\alpha^{n_\alpha}\bigl((\sigma^\ell,\rho^\ell)_{\ell\leq n}\bigr) \Bigr\ra,
\label{ch45lhscor}
\end{equation}
where $W_\alpha$'s were defined in (\ref{ch45WA}). Let us now compute the right hand side of (\ref{ch45main}). By (\ref{fgzeros}), the coordinates $(\sigma^{\ell},\rho^\ell)$ for $\ell>n$ are not present in all the functions $D_\ell$ defined above and we will continue to write them as functions of $(\sigma,\rho)$ and $(\sigma^\ell,\rho^\ell)_{\ell\leq n}$ only. Then, it is easy to see that the denominator on the right hand side of (\ref{ch45main}) is equal to $\la\exp D(\sigma,\rho)\ra_{\hspace{-0.3mm}\mathunderscore}^m$ and the sum in the exponent in the numerator equals $\sum_{\ell=1}^{m} D_{\ell}(\sigma^{\ell},\rho^\ell)$, where $D$ and $D_\ell$ depend implicitly on $(\sigma^\ell,\rho^\ell)_{\ell\leq n}$ and are defined exactly as above Theorem \ref{ch45Th1}.

Since the function $\Phi$ and denominator do not depend on  $(\sigma^{\ell},\rho^\ell)_{\ell>n}$, integrating the numerator in the coordinate $(\sigma^{\ell},\rho^\ell)$ for $\ell\in S_\alpha$ produces a factor 
$
\la I_{B_\alpha}(\sigma,\rho) \exp D(\sigma,\rho)  \ra_{\hspace{-0.3mm}\mathunderscore}.
$
For each $\alpha\in \A$, we have $|S_\alpha| = n_\alpha$ such coordinates and, therefore, the right hand side of (\ref{ch45main}) 
is equal to
$$
\e\Bigl\la
\frac{\Phi(R^n) \exp \sum_{\ell=1}^{n} D_{\ell}(\sigma^{\ell},\rho^\ell)}
{\la\exp D(\sigma,\rho)\ra_{\hspace{-0.3mm}\mathunderscore}^n}
\! \prod_{\alpha\in\A} \hspace{-0.3mm} \Bigl(\hspace{-0.1mm}
\frac{\la I_{B_\alpha}(\sigma,\rho) \exp D(\sigma,\rho)\ra_{\hspace{-0.3mm}\mathunderscore}}
{\la \exp D(\sigma,\rho )\ra_{\hspace{-0.3mm}\mathunderscore}}
\Bigr)^{n_\alpha}\Bigr\ra.
$$
Comparing this with (\ref{ch45lhscor}) and recalling the notation (\ref{ch45TA}) proves (\ref{ch45main2}) for
$$
\varphi(R^n, W) =\Phi(R^n) \prod_{\alpha\in \A} W_\alpha^{n_\alpha}.
$$
The general case then follows by approximation. First, we can approximate a continuous function $\phi$ on $[0,1]^{|\A|}$ by polynomials to obtain (\ref{ch45main2}) for products $\Phi(R^n) \phi(W)$. This, of course, implies the result for continuous functions $\varphi(R^n,W)$ and then for arbitrary bounded measurable functions.
\qed

\medskip
\noindent
It will be convenient to rewrite the above invariance properties in the case when the functions $\varphi$ and the partition depend on different number of replicas $\sigma^1,\ldots,\sigma^n$ and $\rho^1,\ldots,\rho^m$. If we suppose that $m\leq n$ and sets the functions  
$
g_{m+1}=\ldots=g_n=0
$
then the invariance property can be rewritten as follow. First, as before we define
$$
F(\tau,\sigma^1,\ldots,\sigma^n) = f_1(\tau\cdot\sigma^1)+\ldots+f_n(\tau\cdot\sigma^n).
$$
Also, as before, for $1\leq \ell\leq n$ we write
$$
F_{\ell}(\tau,\sigma^1,\ldots,\sigma^n) = F(\tau,\sigma^1,\ldots,\sigma^n)
 - f_{\ell}(\tau\cdot\sigma^{\ell})+ \e \bigl \la f_{\ell}\bigl(\sigma^1\cdot \sigma^2 \bigr) \bigr\ra
$$
and, for $\ell\geq n+1$, we write
$$
F_{\ell}(\tau,\sigma^1,\ldots,\sigma^n) = F(\tau,\sigma^1,\ldots,\sigma^n).
$$
Since $m\leq n$, we have
$$
G(\tau,\rho^1,\ldots,\rho^m) = g_1(\tau\cdot\rho^1)+\ldots+g_m(\tau\cdot\rho^m)
$$
and, for $1\leq \ell\leq m$,
$$
G_{\ell}(\tau,\rho^1,\ldots,\rho^m) = G(\tau,\rho^1,\ldots,\rho^m) 
 - g_{\ell}(\tau \cdot\rho^{\ell})+ \e \bigl \la g_{\ell}\bigl(\rho^1\cdot\rho^2 \bigr) \bigr\ra.
$$
For $\ell\geq m+1$, we now have
$$
G_{\ell}(\tau,\rho^1,\ldots,\rho^m)=G(\tau,\rho^1,\ldots,\rho^m).
$$
We will decouple the remaining notation as follows. Let us denote
\begin{align}
D^1(\sigma)
&= 
F(\sigma,\sigma^1,\ldots,\sigma^n) + \kappa G(\sigma,\rho^1,\ldots,\rho^m),
\\
D^2(\rho)
&= 
 \frac{1}{\kappa} F(\rho,\sigma^1,\ldots,\sigma^n) + G(\rho,\rho^1,\ldots,\rho^m)
\end{align}
and, for $\ell\geq 1$, let us denote
\begin{align}
D_\ell^1(\sigma)
& = 
F_\ell(\sigma,\sigma^1,\ldots,\sigma^n) + \kappa G(\sigma,\rho^1,\ldots,\rho^m),
\\
D_\ell^2(\rho)
&= 
 \frac{1}{\kappa} F(\rho,\sigma^1,\ldots,\sigma^n) + G_\ell(\rho,\rho^1,\ldots,\rho^m).
\end{align}
All of these functions now implicitly depend on $\sigma^1,\ldots,\sigma^n$ and $\rho^1,\ldots,\rho^m$. Then
\begin{align}
D(\sigma,\rho) &=D^1(\sigma)+D^2(\rho),
\\
D_\ell(\sigma,\rho) &=D_\ell^1(\sigma)+D_\ell^2(\rho)
\end{align}
for all $\ell\geq 1.$ With this notation,
$$
\bigl\la\exp D(\sigma,\rho)\bigr\ra_{\hspace{-0.3mm}\mathunderscore}
=
\bigl\la\exp D^1(\sigma)\bigr\ra_{\hspace{-0.3mm}\mathunderscore}
\bigl\la\exp D^2(\rho)\bigr\ra_{\hspace{-0.3mm}\mathunderscore}
$$
and
$$
\sum_{\ell=1}^n D_\ell(\sigma^\ell,\rho^\ell)
=
\sum_{\ell=1}^n D_\ell^1(\sigma^\ell)+\sum_{\ell=1}^n D_\ell^2(\rho^\ell).
$$
When $\varphi(R^n,W)$ and the partition $(B_\alpha)_{\alpha\in\A}$ do not depend on the coordinates $\rho^{m+1},\ldots,\rho^n$, the factors $\exp D_\ell^2(\rho^\ell)$ for $\ell\geq m+1$ in the numerator in (\ref{ch45main2}) can be integrated with respect to $G_2$ and cancelled out with the corresponding factors $\la\exp D^2(\rho)\ra_{\hspace{-0.3mm}\mathunderscore}$ in the denominator. Therefore, Theorem \ref{ch45Th2} can be rewritten as follow.
\begin{theorem}\label{ch45Th3}
Let $\Phi$ be a bounded measurable function of the overlaps $R$ of $\sigma^1,\ldots,\sigma^n$ and $\rho^1,\ldots,\rho^m$ and the weights $W$ in (\ref{ch45WA}), which are also defined in terms of the partition that depends only on these replicas. Then,
\begin{equation}
\e \bigl\la  \varphi(R,W) \bigr\ra
=
\e\biggl\la
\frac{ \varphi(R, T(W)) \exp \bigl(\sum_{\ell=1}^n D_\ell^1(\sigma^\ell)+\sum_{\ell=1}^m D_\ell^2(\rho^\ell) \bigr)}
{\bigl\la\exp D^1(\sigma)\bigr\ra_{\hspace{-0.3mm}\mathunderscore}^n \bigl\la\exp D^2(\rho)\bigr\ra_{\hspace{-0.3mm}\mathunderscore}^m}
\biggr\ra.
\label{ch45main3}
\end{equation}
where the averages $\la\,\cdot\,\ra_{\hspace{-0.3mm}\mathunderscore}$ in the denominator are with respect to $G_1$ or $G_2$.
\end{theorem}

\section{Joint clustering for large values of overlaps}

In this section, we will give first application of the invariance properties above and show that the overlaps within systems and between the two systems satisfy a joint clustering property for large values of the overlaps.

First, let us mention one standard property of the generic models, namely, that the distributions of the absolute values of the overlaps $|\tilde{\sigma}^1\cdot\tilde{\sigma}^2|$ and $|\tilde{\rho}^1\cdot\tilde{\rho}^2|$ within the two system converge weakly to some measures $\mu_1$ and $\mu_2$ on $[0,1]$, called the Parisi measures. This is a standard consequence of the Parisi formula for the free energy (see e.g. Talagrand \cite{PM},  Theorem 14.11.6 in \cite{SG2-2} or Section 3.7 in \cite{SKmodel}). These measures also appear as unique minimizers in the Parisi formula, which will come up in Section \ref{SecCoupled} below. From now on we will denote by
\begin{equation}
c_1 = \inf \mathrm{supp}\, \mu_1,\,\,
c_2 = \inf \mathrm{supp}\, \mu_2
\label{smallestcs}
\end{equation}
the smallest points in the support of the Parisi measures $\mu_1$ and $\mu_2$.

Let us consider $n$ replicas $\sigma^1,\ldots,\sigma^n$ from $G_1$ and $m$ replicas $\rho^1,\ldots,\rho^m$ from $G_2$. To simplify notation, let us denote them by
$$
\tau^\ell = \sigma^\ell \mbox{ for } 1\leq \ell\leq n,\,\,
\tau^{n+\ell} = \rho^\ell \mbox{ for } 1\leq \ell\leq m. 
$$
We will prove that any of these points can be `duplicated' in a certain sense that will be explained below (see the first remark below Theorem \ref{ch45ThObs}) and, for certainty, we will fix that point to be $\sigma^1$. Then we will denote the duplicate point $\sigma^{n+1}$ by $\tau^{n+m+1}$, and $\tau^{n+m+1}$ will represent $\sigma^{n+1}.$

Let us consider an overlap array
$$
R^{n+m} = \bigl(\tau^\ell\cdot\tau^{\ell'}\bigr)_{1\leq \ell<\ell'\leq n+m}
$$
and an array of some fixed parameters
\begin{equation}
A=\bigl(a_{\ell,\ell'}\bigr)_{1\leq \ell<\ell'\leq n+m}.
\end{equation}
Given $\eps>0$, we will write $x\approx a$ to denote that $a-\eps<x<a+\eps$ and $R^{n+m}\approx A$ to denote that $R_{\ell,\ell'} \approx a_{\ell,\ell'}$ for all $1\leq \ell< \ell' \leq n+m$ and, for simplicity of notation, we will keep the dependence of $\approx$ on $\eps$ implicit. Below, the matrix $A$ will be used to describe a set of constraints such that the overlaps in $R^{n+m}$ can take values close to $A$, 
\begin{equation}
\e\bigl\la I\bigl(R^{n+m} \approx A\bigr)\bigr\ra >0,
\label{ch45support}
\end{equation}
for a given $\eps>0$. Let us consider the quantity 
\begin{equation}
a^* = \max \bigl(|a_{1,2}|,\ldots, |a_{1,n+m}|\bigr).
\label{astar1}
\end{equation}
Recall that $c_1$ is the smallest point in the support of the Parisi measure $\mu_1$ defined in (\ref{smallestcs}) and consider any
\begin{equation}
x\geq \max(c_1, a^*).
\label{condonx}
\end{equation}
Then the following duplication property holds.
\begin{theorem}[Duplication I] \label{ch45ThObs} 
Given $\eps>0$, if (\ref{ch45support}) and (\ref{condonx}) hold then
\begin{equation}
\e\Bigl\la
I \Bigl(
R^{n+m}\approx A, \tau^{\ell}\cdot\tau^{n+m+1} \approx a_{1,\ell} \mbox{ for } 2\leq \ell\leq n+m, \bigl|\tau^1\cdot \tau^{n+m+1} \bigr| < x +\eps
\Bigr)
\Bigr\ra
>0.
\label{ch45extend}
\end{equation}
\end{theorem}
\noindent
In other words, if replicas $\sigma^1,\ldots,\sigma^n$ and $\rho^1,\ldots,\rho^m$ form some admissible configuration then, with positive probability, we can find a configuration with an additional point $\tau^{n+m+1}$, in this case $\sigma^{n+1}$, which has (approximately) the same overlap as $\sigma^1$ with all other replicas and at the same time its overlap with $\sigma^1$ is not too big, $|\sigma^1\cdot \sigma^{n+1}| < x +\eps$.

\medskip
\noindent
\textbf{Remark.}
This result will be used in the following way. Suppose that the array $A$ is in the support of the distribution of $R^{n+m}$ under $\e(G_1\times G_2)^{\otimes \infty}$, which means that (\ref{ch45support}) holds for all $\eps>0$. Then (\ref{ch45extend}) also holds for all $\eps>0$. This means that the support of the distribution of $R^{n+m+1}$ with $\tau^{n+m+1}=\sigma^{n+1}$ under $\e(G_1\times G_2)^{\otimes \infty}$ intersects the event in (\ref{ch45extend}) for every $\eps>0$ and, hence, it contains a point in the set
\begin{equation}
\bigl\{
R^{n+m}= A, \tau^{\ell}\cdot\sigma^{n+1} = a_{1,\ell} \mbox{ for } 2\leq \ell\leq n+m, \bigl|\sigma^1\cdot \sigma^{n+1} \bigr| \leq x
\bigr\},
\label{ch45Aplus}
\end{equation}
since the support is compact. Often when we say below that a point can be duplicated, it does not mean that we keep the same points and add another one, but that if a certain configuration is admissible (in the support of the overlaps) then a duplicated configuration is admissible and one can find possibly different points (or even for a different realization of the measures $G_1$ and $G_2$) with such duplicated overlaps.

\medskip
\noindent
\textbf{Remark.} Theorem \ref{ch45ThObs} also holds if we would like to duplicate one of the points $\rho^1,\ldots,\rho^m$, let us say $\rho^1$, but in this case we have to replace the condition (\ref{condonx}) by
\begin{equation}
x\geq \max(c_2, b^*),
\label{condonx-second}
\end{equation}
where $c_2$ is the smallest point in the support of the Parisi measure $\zeta_2$ and $b^* = \max_{\ell\not = n+1} |a_{\ell,n+1}|$. 

\medskip
\noindent
We will start with the following simple result. 
\begin{lemma}\label{LemDS2}
If $\mu_1(A)>0$ then with probability one for $G_1$-almost all $\sigma^1,$ $G_1(\sigma^2: |\sigma^1\cdot \sigma^2| \in A)>0$.
\end{lemma}
\noindent
Of course, the same statement holds for the measure $G_2$ under the assumption $\mu_2(A)>0$.

\medskip
\noindent
\textbf{Proof.} We have $a=\mu_1(A^c)<1.$  First of all, using the Ghirlanda-Guerra identities, \setstretch{1}
\noindent
\begin{align*}
\e \bigl\la I\bigl( |\sigma_{1}\cdot\sigma_{\ell}| \in A^c, 2\leq \ell\leq n+1\bigr) \bigr\ra
& =\,
\e \bigl\la I\bigl(|\sigma_{1}\cdot\sigma_{\ell}| \in A^c, 2\leq \ell\leq n\bigr) I\bigl(|\sigma_{1}\cdot\sigma_{n+1}| \in A^c\bigr) \bigr\ra
\\
&=\,
\frac{n-1+a}{n}\smsp \e \bigl \la I\bigl(|\sigma_{1}\cdot\sigma_{\ell}| \in A^c, 2\leq \ell\leq n\bigr) \bigr\ra,
\end{align*}
\setstretch{1.08}
\hspace{-1mm}where $\la\,\cdot\,\ra$ is now the average with respect to $G_1^{\otimes\infty}$. Repeating the same computation, one can show by induction on $n$ that this equals
$$
\frac{(n-1+a)\cdots(1+a) a}{n!}
=
\frac{a(1+a)}{n}\Bigl(1+\frac{a}{2}\Bigr)\cdots\Bigl(1+\frac{a}{n-1}\Bigr).
$$
Using the inequality $1+x\leq e^x$, it is now easy to see that
$$
\e \bigl\la I\bigl(|\sigma_{1}\cdot\sigma_{\ell}| \in A^c, 2\leq \ell\leq n+1\bigr) \bigr\ra
\leq
\frac{a(1+a)}{n} e^{a\log n}
=
\frac{a(1+a)}{n^{1-a}}.
$$
If we rewrite the left hand side using Fubini's theorem then, since $a<1,$ letting $n\to \infty$ implies that 
$$
\lim_{n\to \infty} \e \int\! G_1( \sigma^2: |\sigma^1\cdot \sigma^2| \in A^c)^n \smsp dG_1(\sigma^1) =0.
$$   
This leads to contradiction if we assume that $G_1(\sigma^2: |\sigma^1\cdot \sigma^2| \in A^c)=1$ with positive probability over the choice of $G_1$ and the choice of $\sigma^1$, which finishes the proof.
\qed

\medskip
\noindent
\textbf{Proof of Theorem \ref{ch45ThObs}.}
We will prove (\ref{ch45extend}) by contradiction, so suppose that the left  hand side is equal to zero. We will apply Theorem \ref{ch45Th2} with ${\cal A}=\{1,2\}$ and the partition of $H^2$,
$$
B_1 = \bigl\{(\sigma,\rho) : |\sigma\cdot \sigma^1| \geq x +\eps \bigr\},\, 
B_2 = B_1^c.
$$ 
By (\ref{condonx}), $\mu_1([0,x +\eps))>0$ and, by Lemma \ref{LemDS2}, the weight 
$$
W_2=(G_1\times G_2)(B_2) = G_1(\sigma \,:\, |\sigma\cdot \sigma^1| < x +\eps)>0
$$ 
with probability one. Since $W_1=1-W_2$, we can find $p<1$ and small $\delta>0$ such that 
\begin{equation}
\delta\leq 
\e\Bigl\la
I\bigl(
R^{n+m}\approx A, W_1 < p
\bigr)
\Bigr\ra.
\label{ch45littlec}
\end{equation}
Let us apply Theorem \ref{ch45Th3} with the above partition, the choice of 
\begin{equation}
\varphi(R,W) = I\bigl(R^{n+m}\approx A, W_1 <p \bigr),
\end{equation}
and the choice of functions $f_{1}(s) = t I(|s|\geq x+\eps)$ for $t\geq 0$ and all other functions $f_j$ and $g_j$ equal to zero. Since on the event $\{R^{n+m}\approx A\}$ the overlaps $|\sigma^1\cdot\tau^\ell|< |a_{1,\ell}|+\eps \leq x+\eps$ for all $\ell\geq 2$, the sum in the numerator on the right hand side of (\ref{ch45main3}) will become 
\begin{align}
\sum_{\ell=1}^n D_\ell^1(\sigma^\ell)+\sum_{\ell=1}^m D_\ell^2(\rho^\ell)
& = 
\sum_{\ell=2}^{n} t  I\bigl(|\sigma^1\cdot\sigma^\ell | \geq x+\eps\bigr)
+t  \e \bigl\la I\bigl(|\sigma^1\cdot\sigma^2| \geq x + \eps\bigr)\bigr \ra
\nonumber
\\&
+\sum_{\ell=1}^{m} \frac{t}{k}  I\bigl(|\sigma^1\cdot\rho^\ell | \geq x+\eps\bigr)
\nonumber
= 
t \e \bigl\la I\bigl(|\sigma^1\cdot\sigma^2 |\geq x + \eps\bigr)\bigr \ra =: t\gamma.
\nonumber
\end{align}
Since the denominator on the right hand side of (\ref{ch45main3}) is greater or equal to $1$, because $D^1,D^2\geq 0$, the equations (\ref{ch45main3}) and (\ref{ch45littlec}) imply
\begin{equation}
\delta 
\leq
\e\Bigl\la
I\bigl(R^{n+m}\approx A, (T_t(W))_1< p \bigr)  \smsp e^{ t\gamma} 
\Bigr\ra.
\label{ch45littlec2}
\end{equation}
Recalling the definition of the map $T(W)$ in (\ref{ch45TA}), our choice of $B_1$ and $f_1$ implies that
\begin{equation}
(T_t(W))_1 = \frac{W_1 e^t}{W_1e^t + 1-W_1}.
\label{ch45TtW}
\end{equation}
In the average $\la\,\cdot\,\ra$ on the right hand side of (\ref{ch45littlec2}) let us fix $\tau^2,\ldots, \tau^{n+m}$ and consider the average with respect to $\sigma^1$ first. Clearly, on the event $\{R^{n+m}\approx A\}$ such average will be taken over the set
\begin{equation}
\Omega(\tau^2,\ldots,\tau^{n+m}) = \bigl\{\sigma: \sigma\cdot \tau^{\ell} \approx a_{1,\ell} \mbox{ for }2\leq \ell\leq n+m\bigr\}.
\label{ch45Omega}
\end{equation}
Suppose that with positive probability over the choice of the measure $G_1\times G_2$ and replicas  $\tau^2,\ldots,$ $\tau^{n+m}$ satisfying the constraints in $A$, i.e. $\tau^{\ell}\cdot \tau^{\ell'}\approx a_{\ell,\ell'}$ for $2\leq \ell,\ell'\leq n+m$,  we can find two points $\sigma'$ and $\sigma''$ in the support of $G_1$ that belong to the set $\Omega(\tau^2,\ldots,\tau^{n+m})$ and such that  $|\sigma'\cdot \sigma'' |< x + \eps.$ This would then imply
\begin{equation}
\e\Bigl\la
I \Bigl(
R^{n+m}\approx A, \tau^{\ell}\cdot\sigma^{n+1} \approx a_{1,\ell} \mbox{ for } 2\leq \ell\leq n+m, \bigl|\sigma^1\cdot \sigma^{n+1} \bigr| < x +\eps
\Bigr)
\Bigr\ra
>0,
\label{ch45extendAG}
\end{equation}
because for $(\sigma^1,\sigma^{n+1})$ in a small neighborhood of $(\sigma',\sigma'')$ the vector 
$(\tau^1,\ldots,\tau^{n+m},\sigma^{n+1})$ would belong to the event on the left hand side,
$$
\bigl\{
R^{n+m}\approx A, \tau^{\ell}\cdot\sigma^{n+1} \approx a_{1,\ell} \mbox{ for } 2\leq \ell\leq n+m, \bigl|\sigma^1\cdot \sigma^{n+1} \bigr| < x +\eps
\bigr\}.
$$
Since we assumed that the left hand side of (\ref{ch45extendAG}) is equal to zero, we must have that, for almost all choices of the measure $G_1\times G_2$ and replicas $\tau^2,\ldots, \tau^{n+m}$ satisfying the constraints in $A$, any two points $\sigma',\sigma''$ in the support of $G_1$ that belong to the set $\Omega(\tau^2,\ldots,\tau^{n+m})$ satisfy $|\sigma'\cdot \sigma''| \geq x + \eps.$ In other words, given a point $\sigma'$, we can not find a point $\sigma''$ in the support of $G_1$ such that  $|\sigma'\cdot \sigma''| < x + \eps.$

Let us also recall that in (\ref{ch45littlec2}) we are averaging over $\sigma^1$ that satisfy the condition $(T_t(W))_1 < p.$ This means that if we fix any such $\sigma'$ in the support of $G_1$ that satisfies this condition and belongs to the set (\ref{ch45Omega}) then the Gibbs average in $\sigma^1$ will be taken over the set
$$
B_1 = B_1(\sigma') = \bigl\{\sigma'': |\sigma'\cdot \sigma''| \geq x +\eps \bigr\}
$$
of measure $W_1 = W_1(\sigma') = G_1(B_1(\sigma'))$ that satisfies $(T_t(W))_1<p.$ It is easy to check that the inequality 
$$
(T_t(W))_1 = \frac{W_1 e^t}{W_1e^t + 1-W_1} < p
$$ 
implies that $W_1 \leq (1-p)^{-1}e^{-t}$. This means that the average on the right hand side of (\ref{ch45littlec2}) over $\sigma^1$ for fixed $\tau^2,\ldots, \tau^{n+m}$ is bounded by $(1-p)^{-1} e^{-t} e^{t\gamma }$ and, thus, for $t\geq 0$,
$$
0<\delta  
\leq 
\e\Bigl\la
I\bigl (R^{n+m}\approx A, (T_t(W))_1<p \bigr)  e^{ t\gamma} 
\Bigr\ra 
\leq 
(1-p)^{-1} e^{-t(1-\gamma)}.
$$
Since $x\geq c_1$ by the assumption (\ref{condonx}), 
$$
1-\gamma = \e\la I(|\sigma^1\cdot\sigma^2|< x + \eps)\ra 
=\mu_1\bigl([0,x+\eps)\bigr)
> 0,
$$ 
and letting $t\to+\infty$ we arrive at contradiction.
\qed

\medskip
\noindent
Using the above duplication property, we will now prove a joint clustering property for the absolute values of the overlaps from the measures $G_1$ and $G_2$. Let us consider
\begin{equation}
q\geq \max(c_1,c_2).
\label{condonx3}
\end{equation}
Then the following holds.
\begin{theorem}[Clustering]\label{ThCluster}
For $q$ that satisfies (\ref{condonx3}), we have
\begin{equation}
\bigl\{|\tau^1\cdot\tau^2|\geq q\bigr\}
\cap \bigl\{|\tau^1\cdot\tau^3|\geq q\bigr\}
\subseteq
\bigl\{|\tau^2\cdot\tau^3|\geq q\bigr\}
\label{ch44ultra}
\end{equation}
with probability one over choice of any three replicas $\tau^1,\tau^2,\tau^2$ (which could be $\sigma$'s from $G_1$ or $\rho$'s from $G_2$).
\end{theorem}
\noindent
\textbf{Proof.}
The proof is by contradiction. Suppose that (\ref{ch44ultra}) is violated, in which case there exist 
$$
\max(c_1,c_2)\leq a<b\leq c
$$ 
such that the vector $(c,b,a)$ is in the support of $(|\tau^1\cdot\tau^2|, |\tau^1\cdot\tau^3|, |\tau^2\cdot\tau^3|).$ Then, there exists a particular choice of $\eps_1,\eps_2,\eps_3\in\{-1,+1\}$ such that $(\eps_1 c,\eps_2 b,\eps_3 a)$ is in the support of the array of overlaps $(\tau^1\cdot\tau^2, \tau^1\cdot\tau^3, \tau^2\cdot\tau^3).$

Using Theorem \ref{ch45ThObs} repeatedly with the choice of $x=c$, we can duplicate points $\tau^1$ and $\tau^2$ (in the sense described in the remark below Theorem \ref{ch45ThObs}) as many times as we like  by preserving the scalar products with other points and at the same time making sure that no two points have scalar product in absolute value exceeding $c$. This means that, for any $n\geq 1$, there exist points $\tau^3$, $(\tau_\ell^1)_{\ell\leq n}$ and $(\tau_\ell^2)_{\ell\leq n}$ in our Hilbert space $H$ such that
$$
\tau_\ell^1\cdot\tau^3 = \eps_2 b,\,\, 
\tau_\ell^2\cdot\tau^3 = \eps_3 a,\,\, 
\tau_\ell^1\cdot\tau_{\ell'}^2 = \eps_1 c,\,\, 
|\tau_\ell^1\cdot\tau_{\ell'}^1|\leq c,\,\,
 |\tau_\ell^2\cdot\tau_{\ell'}^2|\leq c
$$
for all $\ell,\ell' \leq n$. Let us consider the barycenters of these sets of duplicate points,
$$
\xoverline{\tau}^1 = \frac{1}{n}\sum_{\ell=1}^n \tau_\ell^1,\,\,
\xoverline{\tau}^2 = \frac{1}{n}\sum_{\ell=1}^n \tau_\ell^2.
$$
Then $\xoverline{\tau}^1\cdot \xoverline{\tau}^2 = \eps_1 c,$ $\xoverline{\tau}^1\cdot {\tau}^3 = \eps_2 b$  and $\xoverline{\tau}^2\cdot {\tau}^3 = \eps_3 a$ and
$$
\|\xoverline{\tau}^j\|^2 = \frac{1}{n^2}\sum_{\ell\leq n} \|\tau_{\ell}^j\|^2 + \frac{1}{n^2}\sum_{\ell\not = \ell'} \tau_\ell^j\cdot\tau_{\ell'}^j \leq \frac{n + n(n-1) c}{n^2}.
$$
Therefore, we can write
$$
\|\xoverline{\tau}^1- \eps_1 \xoverline{\tau}^2\|^2 = \|\xoverline{\tau}^1\|^2 + \|\xoverline{\tau}^2\|^2 
- 2\eps_1 \xoverline{\tau}^1\cdot \xoverline{\tau}^2 \leq \frac{2(1 -c)}{n},
$$
which implies that
$
|\eps_2 b- \eps_1 \eps_3  a| 
= |{\tau}^3 \cdot \bigl( \xoverline{\tau}^1 -\eps_1 \xoverline{\tau}^2\bigr)| 
\leq 2n^{-1/2}.
$
Letting $n\to\infty$ yields $\eps_2 b = \eps_1 \eps_3  a$, which contradicts the assumption $0\leq a<b$.
\qed

\section{Large values of the overlap}

Let $\mu_1$ and $\mu_2$ be the Parisi measures and let us define
\begin{equation}
q_0 = \inf\bigl\{t \,:\, \beta_1 \mu_1\bigl([0,t)\bigr) \not= \beta_2\mu_2\bigl([0,t)\bigr)\bigr\}.
\label{qnot}
\end{equation}
From now on, let us for certainty suppose that 
\begin{equation}
c_1\leq c_2,
\label{c1not2}
\end{equation}
where $c_1$ and $c_2$ are the smallest points in the support of $\mu_1$ and $\mu_2$ defined in (\ref{smallestcs}). Let us recall that, by Theorem 4 in Chen \cite{ChenChaos}, if $\e(h^j)^2=0$ then $c_j=0.$ On the other hand, by Theorem 14.12.1 in Talagrand \cite{SG2-2}, if $\e(h^j)^2> 0$ then $c_j>0$. We will now show that the overlap between two systems can not take large values.
\begin{theorem}\label{Thq0}
If (\ref{qnot}) and (\ref{c1not2}) hold then
\begin{equation}
\e\bigl\la I\bigl(|\sigma^1\cdot\rho^1| > \max(c_2,q_0) \bigr)\bigr\ra =0.
\label{qnoteq}
\end{equation}
\end{theorem}
\noindent
\textbf{Proof.} Suppose that (\ref{qnoteq}) is violated,
$$
\e\bigl\la I\bigl(|\sigma^1\cdot\rho^1| > \max(c_2,q_0) \bigr)\bigr\ra > 0.
$$
Then, by the definition (\ref{qnot}), we can find $q>q_0$ such that $q\geq c_2$ and such that
$$
\e\bigl\la I\bigl(|\sigma^1\cdot\rho^1| \geq q \bigr)\bigr\ra >0
\,\,\mbox{ and }\,\,
\beta_1 \mu_1([0,q)) \not= \beta_2\mu_2([0,q)).
$$
Indeed, if $q_0\geq c_2$ then we can find such $q$ right above $q_0$, otherwise,  if $q_0<c_2$ then we can take $q=c_2$ since, in this case, $\mu_2([0,c_2))=0$ and $\mu_1([0,c_2))>0$. By Theorem \ref{ThCluster}, on the event $\{|\sigma^1\cdot\rho^1| \geq q\}$, the following equalities in the support of $G_1$ and $G_2$ hold,
\begin{align}
B &= \bigl\{\sigma\,:\, |\sigma \cdot\sigma^1| \geq q\bigr\} =
\bigl\{\sigma\,:\, |\sigma \cdot\rho^1| \geq q\bigr\}, 
\nonumber
\\
B' &= \bigl\{\rho\,:\, |\rho \cdot\sigma^1| \geq q\bigr\} =
\bigl\{\rho\,:\, |\rho \cdot\rho^1| \geq q\bigr\}. 
\label{Usymmetry}
\end{align}
This is the symmetry property that was mentioned in the introduction. Let us denote $W_1 = G_1(B)$ and $W_2 = G_2(B')$. Let us use Theorem \ref{ch45Th1} for $n=1$ and 
$$
\Phi = I\bigl(|\sigma^1\cdot\rho^1| \geq q \bigr).
$$
First, we apply Theorem \ref{ch45Th1} with $g_1=0$ and $f_1(x) = tI(|x|\geq q)$. Since $f_1(\sigma^1\cdot\rho^1) = t$ on the event $\{|\sigma^1\cdot\rho^1| \geq q\}$, we get
$$
0<\e\bigl\la I\bigl(|\sigma^1\cdot\rho^1| \geq q \bigr)\bigr\ra 
=
\e\biggl\la I\bigl(|\sigma^1\cdot\rho^1| \geq q \bigr)
\frac{\exp \bigl(t \e \la I(|\sigma^1\cdot\sigma^2|\geq q)\ra + t/\kappa \bigr)}
{\bigl(W_1 e^t + 1-W_1\bigr)\bigl(W_2 e^{t/\kappa} + 1-W_2\bigr)}
\biggr\ra.
$$
Next, we apply Theorem \ref{ch45Th1} with $f_1=0$ and $g_1(x) = sI(|x|\geq q)$. In this case, we get
$$
0<\e\bigl\la I\bigl(|\sigma^1\cdot\rho^1| \geq q \bigr)\bigr\ra 
=
\e\biggl\la I\bigl(|\sigma^1\cdot\rho^1| \geq q \bigr)
\frac{\exp \bigl( s \e \la I(|\rho^1\cdot \rho^2|\geq q)\ra + \kappa s\bigr)}
{\bigl(W_1 e^{\kappa s} + 1-W_1\bigr)\bigl(W_2 e^{s} + 1-W_2\bigr)}
\biggr\ra.
$$
If we take $t=\kappa s>0$, the right hand sides above can be equal only if
$$
\kappa s \e \la I(|\sigma^1\cdot\sigma^2|\geq q) + s
=
s \e \la I(|\rho^1\cdot\rho^2|\geq q) + \kappa s
$$
or, equivalently,
$
\e\la I(|\rho^1\cdot\rho^2|< q)=\kappa \e \la I(|\sigma^1\cdot\sigma^2|< q).
$
This can be written as $\beta_1 \mu_1([0,q)) = \beta_2\mu_2([0,q))$, which contradicts our choice of $q$ above. Asymmetry of the invariance property with respect to $\beta_1, \beta_2$ turns out to be incompatible with the symmetry expressed in (\ref{Usymmetry}).
\qed

\section{Intermediate values: uncoupled case}

In this section, we continue to assume (\ref{c1not2}) and we will also assume that
\begin{equation}
q_0 \leq c_2,
\label{q0uncoupled}
\end{equation}
where $q_0$ was defined in (\ref{qnot}). We will call this \emph{uncoupled case}, because this means that either $c_1<c_2$ or if $c_1=c_2$ then the measures $\beta_1\mu_1$ and $\beta_2\mu_2$ are immediately different to the right of $c_1=c_2.$ We will treat the \emph{coupled case} later by very different methods, i.e. when $c_1=c_2$ and the measures $\mu_1$ and $\mu_2$ are equal on some non-trivial interval $[c_1,q_0).$

First of all, if (\ref{q0uncoupled}) holds and if $\e(h^2)^2 = 0$ then, as we mentioned above, $c_2 = 0$ by Theorem 4 in Chen \cite{ChenChaos} and, therefore, Theorem \ref{Thq0} implies that the overlap can not take values other than zero asymptotically. This means that in the rest of this section we can assume that $\e(h^2)^2  > 0$ and $c_2>0$. Moreover, Theorem 14.12.1 in Talagrand \cite{SG2-2} also gives that $c_2>0$ is the smallest point in the support of the distribution of $\rho^1\cdot\rho^2$ and not only $|\rho^1\cdot\rho^2|$, i.e. the overlap is strictly positive for the second system.

\begin{theorem}\label{ThBelowq0}
If (\ref{q0uncoupled}) holds then
\begin{equation}
\e\bigl\la I\bigl(\sigma^1\cdot\rho^1 = \sigma^1\cdot\rho^2 \bigr)\bigr\ra =1.
\label{qnoteqAg}
\end{equation}
\end{theorem}
\noindent
\textbf{Proof.}
Suppose there exist $x$ and $q_1\not = q_2$ such that $(q_1,q_2,x)$ is in the support of
$$
(\sigma^1\cdot\rho^1, \sigma^1\cdot\rho^2,\rho^1\cdot\rho^2).
$$
By Theorem \ref{Thq0}, we must have $|q_1|, |q_2| \leq c_2$. By the comment above, $x\geq c_2>0.$ Using Theorem \ref{ch45ThObs} repeatedly, we can duplicate points $\rho^1$ and $\rho^2$ (in the sense described in the remark below Theorem \ref{ch45ThObs}) as many times as we like. This means that, for any $n\geq 1$, there exist points $\sigma^1$, $(\rho_\ell^1)_{\ell\leq n}$ and $(\rho_\ell^2)_{\ell\leq n}$ in our Hilbert space $H$ such that
$$
\rho_\ell^1\cdot\sigma^1 = q_1, \rho_\ell^2\cdot\sigma^1 = q_2, \rho_\ell^1\cdot\rho_{\ell'}^2 = x, |\rho_\ell^1\cdot\rho_{\ell'}^1|\leq x, |\rho_\ell^2\cdot\rho_{\ell'}^2|\leq x 
$$
for all $\ell,\ell' \leq n$. Let us consider the barycenters of these sets of duplicate points,
$$
\xoverline{\rho}^1 = \frac{1}{n}\sum_{\ell=1}^n \rho_\ell^1,\,\,
\xoverline{\rho}^2 = \frac{1}{n}\sum_{\ell=1}^n \rho_\ell^2.
$$
Then $\xoverline{\rho}^1\cdot \xoverline{\rho}^2 = x,$ $\xoverline{\rho}^1\cdot {\sigma}^1 = q_1$  and $\xoverline{\rho}^2\cdot \sigma^1 = q_2$ and
$$
\|\xoverline{\rho}^j\|^2 = \frac{1}{n^2}\sum_{\ell\leq n} \|\rho_{\ell}^j\|^2 + \frac{1}{n^2}\sum_{\ell\not = \ell'} \rho_\ell^j\cdot\rho_{\ell'}^j \leq \frac{n + n(n-1) x}{n^2}.
$$
Therefore, we can write
$$
\|\xoverline{\rho}^1- \xoverline{\rho}^2\|^2 = \|\xoverline{\rho}^1\|^2 + \|\xoverline{\rho}^2\|^2 
- 2 \xoverline{\rho}^1\cdot \xoverline{\rho}^2 \leq \frac{2(1 -x)}{n},
$$
which implies that
$$
|q_1- q_2| 
= |\sigma^1 \cdot \bigl( \xoverline{\rho}^1 - \xoverline{\rho}^2\bigr)| 
\leq 2n^{-1/2}.
$$
Letting $n\to\infty$ contradicts the assumption that $q_1\not = q_2$.
\qed

\medskip
Next, we will prove another duplication property. As before, let us consider $n$ replicas $\sigma^1,\ldots,\sigma^n$ from $G_1$ and $m$ replicas $\rho^1,\ldots,\rho^m$ from $G_2$. To simplify notation, let us denote them by
$$
\tau^\ell = \sigma^\ell \mbox{ for } 1\leq \ell\leq n,\,\,
\tau^{n+\ell} = \rho^\ell \mbox{ for } 1\leq \ell\leq m. 
$$
We will prove that $\sigma^1$ can be duplicated in the sense very similar to the first duplication property in Theorem \ref{ch45ThObs}, only here the overlaps between $\rho$'s and $\sigma$'s will play a secondary role, due to Theorem \ref{ThBelowq0}. Let us consider an overlap array
$$
R^{n+m} = \bigl(\tau^\ell\cdot\tau^{\ell'}\bigr)_{1\leq \ell<\ell'\leq n+m}
$$
and an array of some fixed parameters
\begin{equation}
A=\bigl(a_{\ell,\ell'}\bigr)_{1\leq \ell<\ell'\leq n+m}.
\end{equation}
As above, given $\eps>0$, we will write $x\approx a$ to denote that $a-\eps<x<a+\eps$ and $R^{n+m}\approx A$ to denote the same approximate equality element-wise. We will assume that 
\begin{equation}
\e\bigl\la I\bigl(R^{n+m} \approx A\bigr)\bigr\ra >0
\label{ch45supportD}
\end{equation}
for all $\eps>0$, i.e. $A$ is in the support of the distribution of $R^{n+m}$. Let us consider the quantity 
\begin{equation}
a^* = \max \bigl(|a_{1,2}|,\ldots, |a_{1,n}|\bigr)
\end{equation}
if $n\geq 2$ and set $a^* = c_1$ if $n=1$. The difference with the duplication property in Theorem \ref{ch45ThObs} above is that now $a^*$ is determined only by the overlaps between $\sigma^1$ and other $\sigma$'s and we ignore the overlaps $a_{1,\ell}$ for $\ell\geq n+1$ between $\sigma^1$ and $\rho$'s. Consider the event
$$
A^+ = \Bigl\{
\sigma^{n+1}\cdot\rho^{\ell} \approx a_{1,n+\ell} \mbox{ for } 1\leq \ell\leq m,  \sigma^{n+1}\cdot\sigma^{\ell} \approx a_{1,\ell} \mbox{ for } 2\leq \ell\leq n, \bigl|\sigma^1\cdot \sigma^{n+1} \bigr| < a^* +\eps
\Bigr\}.
$$
The following duplication property holds.
\begin{theorem}[Duplication II] \label{ch45ThObsD} 
If (\ref{q0uncoupled}) holds and if the array $A$ satisfies (\ref{ch45supportD}) then
\begin{equation}
\e\bigl\la
I \bigl( \bigl\{R^{n}\approx A\bigr\} \cap A^+ \bigr)
\bigr\ra
>0
\label{ch45extendD}
\end{equation}
for all small enough $\eps>0.$
\end{theorem}
\noindent
Again, this can be reinterpreted by saying that if $A$ is in the support of $R^{n+m}$ that the support of the overlaps that include additional replica $\sigma^{n+1}$ contains a points in 
\begin{equation}
\bigl\{
R^{n}= A, \sigma^{n+1}\cdot\rho^{\ell} = a_{1,n+\ell} \mbox{ for } 1\leq \ell\leq m,  \sigma^{n+1}\cdot\sigma^{\ell} = a_{1,\ell} \mbox{ for } 2\leq \ell\leq n, \bigl|\sigma^1\cdot \sigma^{n+1} \bigr| \leq a^*
\bigr\}.
\label{ch45AplusD}
\end{equation}

\medskip
\noindent
\textbf{Proof of Theorem \ref{ch45ThObsD}.}
We will prove (\ref{ch45extendD}) by contradiction, so suppose that the left  hand side is equal to zero. We will apply Theorem \ref{ch45Th2} with ${\cal A}=\{1,2\}$ and the partition of $H^2$,
$$
B_1 = \bigl\{(\sigma,\rho) : |\sigma\cdot \sigma^1| \geq a^* +\eps \bigr\},\, 
B_2 = B_1^c.
$$ 
Since $a^*\geq c_1$, $\mu_1([0,a^* +\eps))>0$ and, by Lemma \ref{LemDS2}, the weight 
$$
W_2=(G_1\times G_2)(B_2) = G_1(\sigma \,:\, |\sigma\cdot \sigma^1| < a^* +\eps)>0
$$ 
with probability one. Therefore, we can find $p<1$ and small $\delta>0$ such that 
\begin{equation}
\delta\leq 
\e\Bigl\la
I\bigl(
R^{n+m}\approx A, W_1 < p
\bigr)
\Bigr\ra.
\label{ch45littlecD}
\end{equation}
Let us apply Theorem \ref{ch45Th3} with the above partition, the choice of 
\begin{equation}
\varphi(R,W) = I\bigl(R^{n+m}\approx A, W_1 <p \bigr),
\end{equation}
and the choice of functions $f_{1}(x) = t I(|x|\geq a^*+\eps)$ for $t\geq 0$ and all other functions $f_j$ and $g_j$ equal to zero. First of all, recalling the definition of the map $T(W)$ in (\ref{ch45TA}), our choice of the set  $B_1$ and functions $f_j$ and $g_j$ implies that
\begin{equation}
(T_t(W))_1 = \frac{W_1 e^t}{W_1e^t + 1-W_1}.
\label{ch45TtWD}
\end{equation}
On the event $R^{n+m}\approx A$, we have $|\sigma^1\cdot \sigma^\ell |< |a_{1,\ell}|+\eps \leq a^* +\eps$ for $2\leq \ell\leq n$ and $f_1(\sigma^1\cdot \sigma^\ell) = 0.$ Therefore, the exponent in the numerator in (\ref{ch45main3}) equals
$$
\exp t \e \bigl\la I\bigl(|\sigma^1\cdot\sigma^2 |\geq a^* + \eps\bigr)\bigr \ra 
\,
\exp \sum_{\ell=1}^m \frac{t}{\kappa} I\bigl(|\sigma^1\cdot\rho^\ell |\geq a^* + \eps\bigr).
$$
By Theorem \ref{ThBelowq0}, we have $\rho^\ell \cdot \sigma^1=\rho^1\cdot \sigma^1$ for all $\ell$, so this is equal to
$$
\exp t \e \bigl\la I\bigl(|\sigma^1\cdot\sigma^2 |\geq a^* + \eps\bigr)\bigr \ra 
\,
\exp m\frac{t}{\kappa} I\bigl(|\sigma^1\cdot\rho^1 |\geq a^* + \eps\bigr).
$$
The denominator in (\ref{ch45main3}) equals 
$$
\bigl\la\exp f_1(\sigma\cdot \sigma^1) \bigr\ra_{\hspace{-0.3mm}\mathunderscore}^n \bigl\la\exp \frac{1}{\kappa} f_1(\rho\cdot\sigma^1)\bigr\ra_{\hspace{-0.3mm}\mathunderscore}^m
=
\bigl(W_1e^t + 1-W_1\bigr)^n \Bigl\la\exp \frac{t}{\kappa} I\bigl(|\sigma^1\cdot\rho |\geq a^* + \eps\bigr)\Bigr\ra_{\hspace{-0.3mm}\mathunderscore}^m.
$$
The behaviour of the second factor in the last two equations will depend on whether
\begin{enumerate}
\item[(i)] $|a_{1,n+1}| \leq a^*$,

\item[(ii)] $|a_{1,n+1}| > a^*$.
\end{enumerate}
In the case (i), on the event $R^{n+m}\approx A$, we have $|\rho^1\cdot \sigma^1|< |a_{1,n+1}|+\eps \leq a^* +\eps$ and, therefore,
$$
\exp \frac{t}{\kappa} I\bigl(|\sigma^1\cdot\rho^1 |\geq a^* + \eps\bigr) = 1.
$$
By Theorem \ref{ThBelowq0}, on the event $R^{n+m}\approx A$, we have $\rho\cdot \sigma^1=\rho^1\cdot \sigma^1$ and, hence, $|\rho\cdot \sigma^1|< a^* +\eps$ for all $\rho$ in the support of $G_2$ and, therefore,
$$
\Bigl\la\exp \frac{t}{\kappa} I\bigl(|\sigma^1\cdot\rho |\geq a^* + \eps\bigr)\Bigr\ra_{\hspace{-0.3mm}\mathunderscore}
=1.
$$
In the case (ii), we can suppose without loss of generality that $\eps>0$ is small enough, so that $|a_{1,n+1}|-\eps> a^*+\eps$. Then, on the event $R^{n+m}\approx A$, $|\rho^1\cdot \sigma^1|\geq |a_{1,n+1}|-\eps > a^* +\eps$ and
$$
\exp m\frac{t}{\kappa} I\bigl(|\sigma^1\cdot\rho^1 |\geq a^* + \eps\bigr) = \exp m\frac{t}{\kappa}.
$$
Again, by Theorem \ref{ThBelowq0}, on the event $R^{n+m}\approx A$, $\rho\cdot \sigma^1=\rho^1\cdot \sigma^1$ and, hence, $|\rho\cdot \sigma^1|\geq a^* +\eps$ for all $\rho$ in the support of $G_2$ and, therefore,
$$
\Bigl\la\exp \frac{t}{\kappa} I\bigl(|\sigma^1\cdot\rho |\geq a^* + \eps\bigr)\Bigr\ra_{\hspace{-0.3mm}\mathunderscore}^m
=\exp m\frac{t}{\kappa}.
$$
In both cases, we see that these terms in the numerator and denominator cancel out. Therefore, if we denote $\gamma = \e \bigl\la I\bigl(|\sigma^1\cdot\sigma^2 |\geq a^* + \eps\bigr)\bigr \ra$, the equations (\ref{ch45main3}) and (\ref{ch45littlecD}) imply
\begin{equation}
\delta 
\leq
\e\Bigl\la
I\bigl(R^{n+m}\approx A, (T_t(W))_1< p \bigr)  \smsp e^{ t\gamma} 
\Bigr\ra.
\label{ch45littlec2D}
\end{equation}
The rest of the argument is identical to the proof of Theorem \ref{ch45ThObsD} and shows that (\ref{ch45littlec2D}) leads to contradiction.
\qed

\medskip
\noindent
Using the above duplication property, we can now prove the following.
\begin{theorem}\label{ThBelowq02}
If (\ref{q0uncoupled}) holds then
\begin{equation}
\e\bigl\la I\bigl(|\sigma^1\cdot\rho^1| = |\sigma^2\cdot\rho^1| \bigr)\bigr\ra =1.
\label{qnoteq2}
\end{equation}
\end{theorem}
\noindent
\textbf{Proof.} If not, then we can find non-negative $x$ and $q_1\not = q_2$ such that $(q_1,q_2,x)$ is in the support of
$$
\bigl(|\sigma^1\cdot\rho^1|, |\sigma^2\cdot\rho^1|, |\sigma^1\cdot \sigma^2|\bigr).
$$
This means that for choice of $\eps_1,\eps_2,\eps_3\in\{-1,+1\}$, $(\eps_1 q_1, \eps_2q_2,\eps_3 x)$ is in the support of
$$
\bigl(\sigma^1\cdot\rho^1, \sigma^2\cdot\rho^1, \sigma^1\cdot \sigma^2\bigr).
$$
Using Theorem \ref{ch45ThObsD} repeatedly, one can show that there exist points $\rho^1$, $(\sigma_\ell^1)_{\ell\leq n}$ and $(\sigma_\ell^2)_{\ell\leq n}$ in our Hilbert space $H$ such that
$$
\sigma_\ell^1\cdot\rho^1 = \eps_1 q_1, \sigma_\ell^2\cdot\rho^1 = \eps_2 q_2, \sigma_\ell^1\cdot\sigma_{\ell'}^2 = \eps_3 x, |\sigma_\ell^1\cdot\sigma_{\ell'}^1|\leq x, |\sigma_\ell^2\cdot\sigma_{\ell'}^2|\leq x 
$$
for all $\ell,\ell' \leq n$. Consider the barycenters,
$$
\xoverline{\sigma}^1 = \frac{1}{n}\sum_{\ell=1}^n \sigma_\ell^1,\,\,
\xoverline{\sigma}^2 = \frac{1}{n}\sum_{\ell=1}^n \sigma_\ell^2.
$$
Then $\xoverline{\sigma}^1\cdot \xoverline{\sigma}^2 = \eps_3 x,$ $\xoverline{\sigma}^1\cdot {\rho}^1 = \eps_1 q_1$  and $\xoverline{\sigma}^2\cdot \rho^1 = \eps_2 q_2$ and
$$
\|\xoverline{\sigma}^j\|^2 = \frac{1}{n^2}\sum_{\ell\leq n} \|\sigma_{\ell}^j\|^2 + \frac{1}{n^2}\sum_{\ell\not = \ell'} \sigma_\ell^j\cdot\sigma_{\ell'}^j \leq \frac{n + n(n-1) x}{n^2}.
$$
Therefore, we can write
$$
\|\xoverline{\sigma}^1- \eps_3 \xoverline{\sigma}^2\|^2 = \|\xoverline{\sigma}^1\|^2 + \|\xoverline{\sigma}^2\|^2 
- 2 \eps_3 \xoverline{\sigma}^1\cdot \xoverline{\sigma}^2 \leq \frac{2(1 -x)}{n},
$$
which implies that
$
|\eps_1 q_1- \eps_2 \eps_3  q_2| 
= |{\rho}^1 \cdot \bigl( \xoverline{\sigma}^1 -\eps_3 \xoverline{\sigma}^2\bigr)| 
\leq 2n^{-1/2}.
$
Letting $n\to\infty$ contradicts the assumption $q_1\not = q_2$.
\qed

\medskip
\noindent
Finally, we get the following.
\begin{theorem}\label{ThBelowq03}
If (\ref{q0uncoupled}) holds then
\begin{equation}
\e\bigl\la I\bigl(|\sigma^1\cdot\rho^1| \geq \sqrt{c_1c_2} \bigr)\bigr\ra =0.
\label{qnoteq3}
\end{equation}
\end{theorem}
\noindent
\textbf{Proof.} Suppose that $q$ is in the support of $|\sigma^1\cdot \rho^1|$. By Theorem \ref{Thq0} and assumption (\ref{q0uncoupled}), $q\leq c_2$. This means that we can apply Theorem \ref{ch45ThObs} to $\rho^1$ with $x=c_2$, since $a^*=q\leq c_2$ in (\ref{astar1}), so we can duplicate $\rho^1$ repeatedly and find $\sigma^1$ and $(\rho^\ell)_{\ell\leq n}$ such that 
$$
|\sigma^1\cdot\rho^\ell |=q, |\rho^\ell\cdot\rho^{\ell'}| \leq c_2
$$
for all $\ell,\ell'\leq n$. Since $c_2$ is the smallest point in the support of $\mu_2$, this means that
$$
|\sigma^1\cdot\rho^\ell |=q, |\rho^\ell\cdot\rho^{\ell'}| = c_2
$$
for all $\ell,\ell'\leq n$.
Again, a more precise statement is that there exist such values in the support of the overlaps of $\sigma^1$ and $(\rho^\ell)_{\ell\leq n}$. Next, we apply Theorem \ref{ch45ThObsD} with $a^*=c_1$ to duplicate $\sigma^1$ repeatedly and find $(\sigma^\ell)_{\ell\leq n}$ and $(\rho^\ell)_{\ell\leq n}$ such that 
$$
|\sigma^\ell\cdot\sigma^{\ell'}| \leq c_1, |\sigma^\ell\cdot\rho^{\ell'} |=q, |\rho^\ell\cdot\rho^{\ell'}| = c_2
$$
for all $\ell,\ell'\leq n$. Since $c_1$ is the smallest point in the support of $\mu_1$, this means that $|\sigma^\ell\cdot\sigma^{\ell'}| = c_1$. By Theorem \ref{ThBelowq0}, there exist $\eps_1,\ldots,\eps_n\in\{-1,+1\}$ such that
$$
|\sigma^\ell\cdot\sigma^{\ell'}| = c_1, \sigma^\ell\cdot\rho^{\ell'}=\eps_\ell q, |\rho^\ell\cdot\rho^{\ell'}| = c_2
$$
for all $\ell,\ell'\leq n$ (i.e. such values are in the support of the corresponding overlaps). Let
$$
\xoverline{\rho} = \frac{1}{n}\sum_{\ell=1}^n \rho^\ell, 
\xoverline{\sigma} = \frac{1}{n}\sum_{\ell=1}^n \eps_\ell \sigma^\ell.
$$
Then we have $\xoverline{\rho} \cdot \xoverline{\sigma} = q$,
$$
\|\xoverline{\sigma}\|^2 \leq \frac{1}{n^2}\sum_{\ell\leq n} \|\sigma_{\ell}\|^2 + \frac{1}{n^2}\sum_{\ell\not = \ell'} |\sigma^\ell\cdot\sigma^{\ell'}| \leq \frac{n + n(n-1) c_1}{n^2}
$$
and
$$
\|\xoverline{\rho}\|^2 \leq \frac{1}{n^2}\sum_{\ell\leq n} \|\rho^{\ell}\|^2 + \frac{1}{n^2}\sum_{\ell\not = \ell'} |\rho^\ell\cdot\rho^{\ell'}| \leq \frac{n + n(n-1) c_2}{n^2}.
$$
Using the Cauchy-Schwarz inequality, $q = \xoverline{\rho} \cdot \xoverline{\sigma} \leq \|\xoverline{\rho}\| \| \xoverline{\sigma}\|$, and using the above bounds and letting $n\to\infty$ shows that $q\leq \sqrt{c_1c_2}.$
\qed

\section{Intermediate values: coupled case}\label{SecCoupled}

In this section, we will consider the complementary case when
\begin{equation}
q_0 > c_2,
\label{q0coupled}
\end{equation}
where $q_0$ was defined in (\ref{qnot}). This means that $c=c_1=c_2$ and the measures $\mu_1$ and $\mu_2$ are equal on some non-trivial interval,
\begin{equation}
\beta_1 \mu_1([0,t]) = \beta_2 \mu_2([0,t])
\mbox{ for all $t\in[c,q_0)$.}
\label{qnotnot}
\end{equation}
We will now go back to the setting of the finite size system on $\{-1,+1\}^N$ and show that the following holds.
\begin{theorem}\label{Thcoupled}
Suppose that (\ref{q0coupled}) holds. Then, for any $\eps>0$ there exists $\delta>0$ such that
\begin{equation}
\lim_{N\to\infty}\e \Bigl\la I\bigl(|\tilde{\sigma}^1\cdot\tilde{\rho}^1| \in [c+\eps,q_0+\delta] \bigr) \Bigr\ra = 0.
\end{equation}
\end{theorem}
Combining this with Theorem \ref{Thq0}, which holds for any subsequential limit, we get the following.
\begin{theorem}\label{ThcoupledFin}
Suppose that (\ref{q0coupled}) holds. Then, for any $\eps>0$,
\begin{equation}
\lim_{N\to\infty}\e \Bigl\la I\bigl(|\tilde{\sigma}^1\cdot\tilde{\rho}^1| \geq c+\eps \bigr) \Bigr\ra = 0.
\end{equation}
\end{theorem}
To prove Theorem \ref{Thcoupled}, we will use Talagrand's analogue of Guerra's replica symmetry breaking bound \cite{Guerra} for coupled system from Theorem 15.7.3 in \cite{SG2-2}. First of all, let us recall the Parisi formula for the free energy for one system,
$$
F_N^j = \frac{1}{N}\e \log \sum_{\sigma} \exp \beta_j \Bigl(H_N(\sigma) + \sum_{i=1}^N h_i^j \sigma_i\Bigr).
$$
Both the Parisi formula and Talagrand-Guerra upper bound can be constructed explicitly, and it is well known that these constructions satisfy certain partial differential equations, as explained, for example, in Section 14.7 in \cite{SG2-2} (see \cite{Jagannath} for a detailed study of the general non-discrete case). Since our arguments below will utilize only the properties expressed by these differential equations (in addition to some well known properties), we will not repeat the explicit constructions here and will only recall the corresponding descriptions in terms of differential equations.

We will abuse notation slightly and write $\mu_j(q) = \mu_j([0,q])$. Let us consider functions $\Phi_j(q,x)$ for $q\in[0,1]$ and $x\in\Reals$ that are solutions of
\begin{equation}
\frac{\partial \Phi_j}{\partial q} = - \frac{\xi''(q)}{2}\Bigl(\frac{\partial^2 \Phi_j}{\partial x^2} + \mu_j(q)\Bigl(\frac{\partial \Phi_j}{\partial x}\Bigr)^2\Bigr)
\label{ParisiPDE}
\end{equation} 
where $\xi(q)$ was defined in (\ref{ch11xidefine}), with the boundary condition at $q=1$ given by
\begin{equation}
\Phi_j(1,x) = \log \ch \bigl(\beta_j(h^j +x)\bigr).
\end{equation}
Let us denote $\theta(q) = q\xi'(q) - \xi(q)$ and define a functional
\begin{equation}
\PP_j(\mu_j) = \e\Phi_j(0,0) - \frac{1}{2}\int_0^1\! \beta_j^2 \mu_j(q)\theta'(q)\,dq,
\end{equation}
where the expectation is in the external field $h^j$. Then the Parisi formula \cite{Parisi79, Parisi} proved by Talagrand in \cite{TPF} (for another proof, see \cite{PPF}) gives that
\begin{equation}
\lim_{N\to\infty} F_N^j = \inf_{\mu_j} \PP(\mu_j).
\end{equation}
It was proved by Auffinger and Chen \cite{AufChen2} that the functional $\PP(\mu_j)$ is strictly convex and the minimizer is unique (see \cite{Jagannath} for another proof). For generic models, this minimizer is precisely the limit of the distribution of the overlaps within systems at the same temperature (see \cite{PM}, Theorem 14.11.6 in \cite{SG2-2} or Section 3.7 in \cite{SKmodel}), so we will continue to denote it by $\mu_j$.

Now, let us consider any $u\in [-1,1]$ and consider the free energy of a coupled system with the overlap constrained to be equal to $u$,
\begin{equation}
F_N(u) = \frac{1}{N}\e \log \sum_{\tilde{\sigma}\cdot\tilde{\rho}=u} \exp \beta_1 \Bigl(H_N(\sigma) + \sum_{i=1}^N h_i^1\sigma_i\Bigr) \exp \beta_2 \Bigl(H_N(\rho) + \sum_{i=1}^N h_i^2 \rho_i\Bigr).
\label{FNU}
\end{equation}
From now on, we will assume that $u\in [0,1]$, because for negative $u$, making the change of variables $\rho\to -\rho$ simply changes $(h_i^1, h_i^2)$ into $(h_i^1, -h_i^2)$, and our arguments will not depend on the choice of the distribution of $(h^1,h^2)$.

We will give an upper bound on (\ref{FNU}) in Proposition \ref{PropT} below, which is just a rephrasing of Proposition 5.1 in Talagrand \cite{TalUltra} with $\lambda=0$ there and $q_{\tau+1}^{1,2} = v\in [0,1].$ We will rewrite this bound for general instead of only discrete parameters using definitions in terms of differential equations, as in (\ref{ParisiPDE}). This last parameter $v$ will represent a value of the overlap up to which the parameters in Talagrand's bound are completely correlated, and after which they are independent. We will take the functions $\xi_{j,j'}$ in that bound to be $\xi_{j,j'}(q) = \beta_j \beta_{j'} \xi(q)$, and we will choose parameters $(n_\ell)$ there in such a way that for the overlaps greater than $v$ they correspond to the c.d.f.s $\mu_1$ and $\mu_2$ of their individual independent systems, and for values of the overlaps less than or equal to $v$ they correspond to some new (improper) c.d.f. $\mu$, which can be arbitrary as long as 
\begin{equation}
\mu(v)\leq \min\bigl(\mu_1(v),\mu_2(v) \bigr).
\label{match}
\end{equation}
Let us consider a function $\Phi_v(q,x)$ for $q\in [0,v]$ and $x\in\Reals$ which is the solution of
\begin{equation}
\frac{\partial \Phi_v}{\partial q} = - \frac{\xi''(q)}{2}\Bigl(\frac{\partial^2 \Phi_v}{\partial x^2} + \mu(q)\Bigl(\frac{\partial \Phi_v}{\partial x}\Bigr)^2\Bigr)
\label{EqPhiv}
\end{equation} 
with the boundary condition at $q=v$ given by
\begin{equation}
\Phi_v(v,x) =\Phi_1(v,x) + \Phi_2(v,x).
\label{Bphiv}
\end{equation}
Let us define a functional
\begin{equation}
\PP(v,\mu) = \e \Phi_v(0,0) - \frac{1}{2}\int_0^{v}\! (\beta_1+\beta_2)^2 \mu(q)\theta'(q)\,dq
- \sum_{j=1}^2 \frac{1}{2}\int_{v}^1\! \beta_j^2 \mu_j(q)\theta'(q)\,dq + \Delta(u,v),
\label{Parisijoint}
\end{equation}
where the expectation is in the external fields $(h^1,h^2)$ and the last term is given by
\begin{equation}
\Delta(u,v) = \beta_1 \beta_2 \bigl( \xi(u) - u \xi'(v) + \theta(v)\bigr).
\end{equation}
Note that for $v=u$, $\Delta(u,v) = 0$. Talagrand's upper bound can be written as follows. 
\begin{proposition}\label{PropT}
For any $u, v\in [0,1]$, if (\ref{match}) is satisfied then
\begin{equation}
F_N(u) \leq \PP(v,\mu).
\label{Tbound}
\end{equation}
\end{proposition}
From now on, we will make the following canonical choice of parameters. We will always take 
\begin{equation}
v< q_0,
\end{equation}
where $q_0$ was defined in (\ref{qnot}), which means that for $q\leq v$ we have
\begin{equation}
\beta_1 \mu_1(q) = \beta_2 \mu_2(q).
\label{equalmus}
\end{equation}
Let us introduce the notation
\begin{equation}
\lambda = \frac{\beta_1}{\beta_1+\beta_2},\,\,
1-\lambda = \frac{\beta_2}{\beta_1+\beta_2} 
\label{deflambda}
\end{equation}
and from now on set the function $\mu(q)$ to be
\begin{equation}
\mu(q): = \lambda \mu_1(q) = (1-\lambda)\mu_2(q)
\label{choicemu}
\end{equation}
for $q\leq v$, which clearly satisfies (\ref{match}). With this choice,
$$
(\beta_1+\beta_2)^2 \mu(q) = (\beta_1+\beta_2) \beta_1 \mu_1(q) = \beta_1^2 \mu_1(q) + \beta_2^2 \mu_2(q)
$$
for $q\leq v<q_0$, where we used (\ref{equalmus}). Therefore, (\ref{Parisijoint}) becomes
\begin{equation}
\PP(v,\mu) = \e\Phi_v(0,0) - \sum_{j=1}^2 \frac{1}{2}\int_{0}^1\! \beta_j^2 \mu_j(q)\theta'(q)\,dq + \Delta(u,v).
\label{Parisijoint2}
\end{equation}
For discrete choices of parameters as in \cite{TalUltra}, one can easily check by looking at the explicit representation of these bounds and using H\"older's inequality that this choice of $\mu$ implies that 
$$
\e\Phi_v(0,0)\leq \e\Phi_1(0,0)+\e \Phi_2(0,0).
$$ 
Our goal will be to show this inequality for arbitrary $\mu_1$ and $\mu_2$ and, moreover, to show that it is strict over a certain range of values of $v$. 
\begin{theorem}\label{ThPhis}
Suppose that (\ref{qnotnot}) holds. Then, for any $\eps>0$ there exists $\delta'>0$ such that
\begin{equation}
\e \Phi_v(0,0)\leq \e \Phi_1(0,0)+ \e \Phi_2(0,0) -\delta'
\label{Phicomp}
\end{equation}
for all $v\in [c+\eps,q_0)$.
\end{theorem}
This immediately implies the following.
\begin{theorem}\label{ThTbst}
Suppose that (\ref{qnotnot}) holds. Then, for any $\eps>0$ there exist $\delta>0$ and $\delta'>0$ such that
\begin{equation}
F_N(u)\leq \sum_{j=1}^2 \PP_j(\mu_j) - \delta'\
\label{FNless}
\end{equation}
for all $u\in [c+\eps,q_0+\delta]$.
\end{theorem}
In a completely standard way, Theorem \ref{Thcoupled} follows from this by classical Gaussian concentration inequalities (see Section 15.7 in \cite{SG2-2}).

\medskip
\noindent
\textbf{Proof of Theorem \ref{ThTbst}.}
For $u\in [c+\eps,q_0)$, let us take $v=u$. Then Proposition \ref{PropT}, (\ref{Parisijoint2}) and (\ref{Phicomp}) imply (\ref{FNless}). For $u\in [q_0,q_0+\delta)$, let us take $v=q_0-\delta$. Then Proposition \ref{PropT}, (\ref{Parisijoint2}) and (\ref{Phicomp}) imply
$$
F_N(u)\leq \sum_{j=1}^2 \PP_j(\mu_j) - \delta' + \Delta(q_0-\delta,q_0+\delta).
$$
Taking $\delta$ small enough we can ensure that $\Delta(q_0-\delta,q_0+\delta)\leq \delta'/2$, so redefining $\delta'$ finishes the proof.
\qed

\medskip
\noindent
The next two results will be proved for a fixed $(h^1,h^2)$. We will begin the proof of Theorem \ref{ThPhis} with the following result. Let us define
\begin{equation}
\Psi_j(q,x) = \frac{1}{\beta_j} \Phi_j(q,x).
\end{equation}
Then the following holds.
\begin{theorem}\label{ThPsiscomp}
If $\beta_1\not =\beta_2$ then, for any fixed $(h^1,h^2)$ and any $q\in [0,1]$, the functions $\Psi_1(q,\,\cdot\,)$ and $\Psi_2(q,\,\cdot\,)$ are not identically equal.
\end{theorem}
\textbf{Proof.} Since $\Phi_j$ satisfy (\ref{ParisiPDE}), $\Psi_j$ satisfy
\begin{equation}
\frac{\partial \Psi_j}{\partial q} = - \frac{\xi''(q)}{2}\Bigl(\frac{\partial^2 \Psi_j}{\partial x^2} + \beta_j \mu_j(q)\Bigl(\frac{\partial \Psi_j}{\partial x}\Bigr)^2\Bigr)
\label{ParisiPDE2}
\end{equation} 
with the boundary condition at $q=1$ given by
\begin{equation}
\Psi_j(1,x) = \frac{1}{\beta_j }\log \ch \bigl(\beta_j(h_j +x)\bigr).
\end{equation}
Let us consider $\psi_j = \frac{\partial \Psi_j}{\partial x}$ and differentiating the above equation in $x$, we see that
\begin{equation}
\frac{\partial \psi_j}{\partial q} = - \frac{\xi''(q)}{2}\Bigl(\frac{\partial^2 \psi_j}{\partial x^2} + 2 \beta_j \mu_j(q)\psi_j \frac{\partial \psi_j}{\partial x}\Bigr)
\label{ParisiPDE3}
\end{equation} 
with the boundary condition at $q=1$ given by
\begin{equation}
\psi_j(1,x) = \mbox{th} \bigl(\beta_j(h_j +x)\bigr).
\end{equation}
It is well known that $|\psi_j(q,x)| \leq 1$ and all partial derivatives of $\psi_j(q,x)$ in $x$ are bounded (see e.g. Propositions 1 and 2 in \cite{AufChen1}). Therefore, we can consider the strong solution of the stochastic differential equation (see e.g. Proposition 8.2.9 in \cite{Karatzas})
\begin{equation}
dX_j(t) = \xi''(t)\beta_j \mu_j(t)\psi_j(t,X_j(t)) dt + \xi''(t)^{1/2} dB(t)
\label{Xjt}
\end{equation}
with $X_j(q) = x$, where $(B(t))_{t\geq 0}$ is a standard Brownian motion. Using It\=o's formula and (\ref{ParisiPDE3}),
\begin{align*}
d\psi_j(t,X_j(t))
= \,& 
\frac{\partial \psi_j}{\partial t} dt + \frac{\partial \psi_j}{\partial x} dX(t)
+ \frac{1}{2} \frac{\partial^2 \psi_j}{\partial x^2} \xi''(t)dt
\\
= \,&
\Bigl(
\frac{\partial \psi_j}{\partial t} +  \xi''(t)\beta_j \mu_j(t)\psi_j \frac{\partial \psi_j}{\partial x}
+ \frac{1}{2} \frac{\partial^2 \psi_j}{\partial x^2} \xi''(t)
\Bigr) dt
+ 
\xi''(t)^{1/2} \frac{\partial \psi_j}{\partial x} dB(t) 
\\
= \,&
\xi''(t)^{1/2} \frac{\partial \psi_j}{\partial x} dB(t).
\end{align*}
Integrating between $q$ and $1$ and taking expectations gives
\begin{equation}
\psi_j(q,x) = \e \psi_j(1,X_j(1)) = \e \mbox{th} \bigl(\beta_j(h_j +X_j(1))\bigr).
\label{Itorepr}
\end{equation}
Let us integrate (\ref{Xjt}) between $q$ and $1$,
$$
X_j(1) - x = \int_{q}^1 \! \xi''(t)\beta_j \mu_j(t)\psi_j(t,X_j(t))\, dt + \int_{q}^1 \xi''(t)^{1/2}\, dB(t).
$$
The first integral is bounded in absolute values by a $\beta_j\xi'(1)$, since $|\psi_j|\leq 1$, and the second term has Gaussian distribution with the variance $\int_q^1 \xi''(t) \,dt = \xi'(1) -\xi'(q)$. Therefore,
$$
\p\Bigl( \bigl| X_j(1) - x\bigr| \geq \gamma \Bigr) \leq e^{-a \gamma^2}
$$
for large $\gamma$ (independent of $x$), where $a$ is some constant that depends on $\beta_1,\beta_2$ and $\xi$. Suppose for certainty that $\beta_1<\beta_2$. First of all, using (\ref{Itorepr}) and the fact that $1-\mathrm{th}(x)$ is decreasing, 
\begin{align*}
1-\psi_1(q,x) 
&= 
1- \e \mbox{th} \bigl(\beta_1(h_1 +X_1(1))\bigr)
\\
&\geq 
\bigl(1- \mbox{th}(\beta_1(h_1 +x +\gamma) \bigr)
\p\bigl( X_1(1) - x \leq \gamma  \bigr)
\\
 & \geq
\bigl(1- \mbox{th}(\beta_1(h_1 +x +\gamma) \bigr)
\bigl( 1-  e^{-a \gamma^2} \bigr).
\end{align*}
Similarly,
\begin{align*}
1-\psi_2(q,x) 
&= 
1- \e \mbox{th} \bigl(\beta_2(h_2 +X_2(1))\bigr)
\\
&\leq 
\bigl(1- \mbox{th}(\beta_2(h_2 +x -\gamma) \bigr)
+
\p\bigl( X_2(1) - x \leq -\gamma  \bigr)
\\
 & \leq
\bigl(1- \mbox{th}(\beta_2(h_2 +x -\gamma) \bigr)
+ e^{-a \gamma^2}.
\end{align*}
Now, let us take $\gamma=\eps x$, where $\eps>0$ is such that $\beta_2(1-\eps)>\beta_1(1+\eps)$. Then, as $x\to+\infty$,
$$
\bigl(1- \mbox{th}(\beta_2(h_2 +x -\gamma) \bigr)
+ e^{-a \gamma^2}
<
\bigl(1- \mbox{th}(\beta_1(h_1 +x +\gamma) \bigr)
\bigl( 1-  e^{-a \gamma^2} \bigr),
$$
since this is equivalent to
$$
\bigl(1- \mbox{th}(\beta_2h_2 +\beta_2(1-\eps)x) \bigr)
+ e^{-a \eps^2x^2}
<
\bigl(1- \mbox{th}(\beta_1h_1 +\beta_1(1+\eps)x) \bigr)
\bigl( 1-  e^{-a \eps^2x^2} \bigr)
$$
and because the two sides have very different asymptotics as $x\to+\infty$ and their ratio goes to zero, using $1-\mbox{th}(x) \sim 2e^{-2x}$. This shows that, for any $q$, the functions $\psi_1(q,\,\cdot\,)$ and $\psi_2(q,\,\cdot\,)$ are not identically equal, which finishes the proof.
\qed

\medskip
\noindent
Recall the definition of $\lambda$ in (\ref{deflambda}) and define
\begin{equation}
\tfi_1(q,x) = \frac{1}{\lambda} \Phi_1(q,x),\,\,
\tfi_2(q,x) = \frac{1}{1-\lambda} \Phi_2(q,x).
\end{equation}
If we recall the definition of $\mu$ in (\ref{choicemu}), one can easily check that 
\begin{equation}
\frac{\partial \tfi_j}{\partial q} = - \frac{\xi''(q)}{2}\Bigl(\frac{\partial^2 \tfi_j}{\partial x^2} + \mu(q)\Bigl(\frac{\partial \tfi_j}{\partial x}\Bigr)^2\Bigr),
\end{equation} 
which is the same equation (\ref{EqPhiv}) satisfied by $\Phi_v(q,x)$ and, moreover, the boundary condition (\ref{Bphiv}) can be rewritten as
\begin{equation}
\Phi_v(v,x) = \lambda \tfi_1(v,x) + (1-\lambda)\tfi_2(v,x).
\label{Bphiv2}
\end{equation}
Theorem \ref{ThPsiscomp} can be expressed by saying that the functions $\tfi_1(q,\,\cdot\,)$ and $\tfi_2(q,\,\cdot\,)$ are not identically equal for any $q$. As a consequence, we will show the following.
\begin{theorem} \label{ThPhiStrict}
If $v\in [c+\eps,q_0)$ then, for any fixed $(h^1,h^2)$ and $q<v$, we have strict inequality
\begin{equation}
\Phi_v(q,x) < \lambda \tfi_1(q,x) + (1-\lambda)\tfi_2(q,x).
\label{Bphiv3}
\end{equation}
\end{theorem}
First, let us show how this implies Theorem \ref{ThPhis}.

\medskip
\noindent
\textbf{Proof of Theorem \ref{ThPhis}.} Let $v_0 := c+\eps$ and consider any $v\in (c+\eps,q_0)$, so that $v_0< v.$ First of all, (\ref{Bphiv3}) implies that 
$$
\Phi_v(v_0,x) < \lambda \tfi_1(v_0,x) + (1-\lambda)\tfi_2(v_0,x) = \Phi_{v_0}(v_0,x).
$$
Both $\Phi_v(q,x)$ and $ \Phi_{v_0}(q,x)$ satisfy the same equation (\ref{EqPhiv}), so $\Phi_v(q,x) \leq \Phi_{v_0}(q,x)$ for all $q\leq v_0$, because this equation preserves monotonicity with respect to the boundary conditions. This is because the case of general $\mu$ can be approximated by step functions (see e.g. Theorem 14.11.2 in \cite{SG2-2}) and, on any interval where $\mu$ is constant, $\exp \mu \Phi_v$ satisfies the heat equation. Using (\ref{Bphiv3}) again implies
$$
\Phi_v(0,0) \leq \Phi_{v_0}(0,0) < \lambda \tfi_1(0,0) + (1-\lambda)\tfi_2(0,0)
=
\Phi_1(0,0) + \Phi_2(0,0).
$$
This holds for any fixed $(h^1,h^2)$ and averaging in $(h^1,h^2)$ yields
$$
\e\Phi_v(0,0) \leq \e \Phi_{v_0}(0,0) < \lambda \e \tfi_1(0,0) + (1-\lambda)\e \tfi_2(0,0)
=
\e \Phi_1(0,0) + \e\Phi_2(0,0).
$$
The strict inequality is uniform over $v\in [c+\eps,q_0)$, which finishes the proof.
\qed

\medskip
It remains to prove Theorem \ref{ThPhiStrict}. Our proof will be based on the variational representation of solutions of the equation (\ref{EqPhiv}) in Theorem 3 in Auffinger, Chen \cite{AufChen2}, which we now recall. 

We will slightly modify their statement to replace the interval $[0,1]$ by $[0,v]$. As before, let $(B(q))_{q\geq 0}$ be a standard Brownian motion and, for $0\leq s< t\leq v$, let $\DD[s,t]$ be the space of all progressively measurable processes $u$ on $[s,t]$ with respect to the filtration generated by $B(q)$ such that $\|u\|_\infty \leq L$ for some arbitrary large enough constant $L$. In \cite{AufChen2}, $L$ was chosen to be $1$ and it could not be smaller that $1$, but it will be convenient not to impose this artificial restriction here. Suppose that $f$ is a solution of (\ref{EqPhiv}), for example, $f=\Phi_v, \tfi_1$ or $\tfi_2$. For any $x\in\Reals$ and $u\in \mathcal{D}[s,t],$ define
\begin{equation}
F^{s,t}(u,x)=\e\bigl( C^{s,t}(u,x)-L^{s,t}(u)\bigr),
\label{defFst}
\end{equation}
where
\begin{align}
C^{s,t}(u,x)
&=
f\Bigl(t,x+\int_s^t\! \xi''(q) \mu(q) u(q)\, dq+\int_s^t\! \xi''(q)^{1/2}\, dB(q)\Bigr),
\nonumber
\\
L^{s,t}(u)
&=
\frac{1}{2}\int_s^t\!  \xi''(q)\mu(q)u(q)^2\, dq.
\end{align} 
With this notation, the following holds (Theorem 3 in \cite{AufChen2}).
\begin{proposition}\label{ThAC}
For any $0\leq s<t\leq v$,
\begin{equation}
f(s,x)=\max_{u\in \DD[s,t]} F^{s,t}(u,x).
\label{ACrepr}
\end{equation}
Moreover, the maximum is attained on 
\begin{equation}
u^*(q)=\frac{\partial f}{\partial x}\bigl(q,X(q) \bigr),
\label{youstar}
\end{equation}
where $(X(q))_{s\leq q\leq t}$ is the strong solution of
\begin{equation}
dX(q) = \xi''(q)\mu(q) \frac{\partial f}{\partial x} \bigl(q,X(q) \bigr) \,dq+\xi''(q)^{1/2} \, dB(q),
\label{ACeq}
\end{equation}
with $X(s)=x.$
\end{proposition}
One can check that $|\frac{\partial f}{\partial x}|$ for $f=\Phi_v, \tfi_1$ or $\tfi_2$ is bounded by $\beta_1+\beta_2,$ so one could in fact take $L=\beta_1+\beta_2$ in the definition of $\DD[s,t]$. 

\medskip
\noindent
\textbf{Proof of Theorem \ref{ThPhiStrict}.} 
First, let us take $t=v$ in Theorem \ref{ThAC} and use it for $f=\Phi_v, \tfi_1$ and $\tfi_2$. Let us denote by $u_v^*, u_1^*$ and $u_2^*$ the corresponding maximizers in (\ref{youstar}). As we mentioned above, Theorem \ref{ThPsiscomp} implies that the functions $\tfi_1(v,\,\cdot\,)$ and $\tfi_2(v,\,\cdot\,)$ are not identically equal. This implies that $u_1^*$ and $u_2^*$ are not identically equal (i.e. they have different trajectories with positive probability). To see this, suppose that they are equal almost surely. Then $X_1$ and $X_2$ are identically equal, $X_1=X_2 = X$, since by (\ref{youstar}) and (\ref{ACeq}),
$$
dX_j(q) = \xi''(q)\mu(q) u_j^*(q) \,dq+\xi''(q)^{1/2} \, dB(q).
$$
Then (\ref{youstar}) would imply that 
$$
\frac{\partial \tfi_1}{\partial x}\bigl(q,X(q) \bigr)= u_1^* =u_2^* =
\frac{\partial \tfi_2}{\partial x}\bigl(q,X(q) \bigr).
$$ 
However, this is impossible because $\frac{\partial \tfi_1}{\partial x}$ and $\frac{\partial \tfi_1}{\partial x}$ are not identically equal for any $q$, while the support of the distribution of $X(q)$ is the whole real line $\Reals$ since, by Girsanov's theorem (Theorem 5.5.1 in \cite{Karatzas}), the distribution of $X(q)$ is Gaussian under some well-defined change of density. 

Let us now show that
\begin{equation}
\Phi_v(s,x) < \lambda \tfi_1(s,x) + (1-\lambda)\tfi_2(s,x).
\label{strictsx}
\end{equation}
Using (\ref{ACrepr}) for $f = \Phi_v$ and (\ref{Bphiv2}), we can write
\begin{align}
\Phi_v(s,x) 
=\,\,\, &
\e\Phi_v\Bigl(v,x+\int_s^v\! \xi''(q) \mu(q) u_v^*(q)\, dq+\int_s^v\! \xi''(q)^{1/2}\, dB(q)\Bigr)
\nonumber
\\
&\, -
\frac{1}{2} \e \int_s^v\!  \xi''(q)\mu(q)u_v^*(q)^2\, dq
\nonumber
\\
=\,\,\,& 
\lambda \Bigl(
\e\tfi_1\Bigl(v,x+\int_s^v\! \xi''(q) \mu(q) u_v^*(q)\, dq+\int_s^v\! \xi''(q)^{1/2}\, dB(q)\Bigr)
\label{step2eq}
\\
&\, -
\frac{1}{2} \e \int_s^v\!  \xi''(q)\mu(q)u_v^*(q)^2\, dq
\Bigr)
\nonumber
\\
& \, +
(1-\lambda) \Bigl(
\e\tfi_2 \Bigl(v,x+\int_s^v\! \xi''(q) \mu(q) u_v^*(q)\, dq+\int_s^v\! \xi''(q)^{1/2}\, dB(q)\Bigr)
\nonumber
\\
&\, -
\frac{1}{2} \e \int_s^v\!  \xi''(q)\mu(q)u_v^*(q)^2\, dq
\Bigr)
\nonumber
\\
\leq \,\,\, &
\lambda \tfi_1(s,x) + (1-\lambda)\tfi_2(s,x),
\nonumber
\end{align}
where the last inequality follows from (\ref{ACrepr}) for $f = \tfi_1$ and $\tfi_2$. Moreover, since we already showed that $u_1^*$ and $u_2^*$ are not identically equal, this inequality will be strict if we can prove that the functional $u\to F^{s,t}(u,x)$ in (\ref{defFst}) is strictly convex. The computation in Proposition 3 in \cite{AufChen2} gives 
$$
\frac{\partial^2}{\partial a^2} F^{s,v}\bigl((1-a)u_1 +a u_2,x\bigr)
\leq
\Bigl(\int_s^v \! \xi''(q)\mu(q)\, dq-1\Bigr)
\e \int_s^v\! \xi''(q)\mu(q)(u_2(q)-u_1(q))^2\, dq.
$$
If $\int_s^v \! \xi''(q)\mu(q)\, dq<1,$ for example, if $|v-s|$ is small, this shows that the functional is strictly convex, so for $s$ close enough to $v$ we obtain (\ref{strictsx}) for all $x$. Because of this, if we take $t=s$ instead of $t=v$ and $s<t$ and repeat the same computation as above, the equality in (\ref{step2eq}) will now become strict inequality and will again yield (\ref{strictsx}). This finishes the proof.
\qed

\section{Small values of the overlap}

Combining Theorem \ref{ThBelowq03} in the uncoupled case with Theorem \ref{ThcoupledFin} in the coupled case gives
\begin{equation}
\lim_{N\to\infty}\e \Bigl\la I\bigl(|\tilde{\sigma}^1\cdot\tilde{\rho}^1| \geq \sqrt{c_1c_2}+\eps \bigr) \Bigr\ra = 0
\label{eqabove}
\end{equation}
for any $\eps>0$. In the coupled case $c_1=c_2=c.$ As we mentioned above, by Theorem 4 in Chen \cite{ChenChaos}, if $\e(h^j)^2=0$ then $c_j=0$, and in this case there is nothing left to prove. If both $\e(h^1)^2 > 0$ and $\e(h^2)^2 > 0$, it remains to appeal to Theorem 7 in Chen \cite{ChenChaos}, which shows that there exists 
$$
\chi \in \bigl[-\sqrt{c_1c_2},\sqrt{c_1c_2} \bigr]
$$ 
that satisfies the following. 
\begin{proposition}
For any $\delta>0$, there exists $\eps>0$ such that
\begin{equation}
\lim_{N\to\infty}\e \Bigl\la I\Bigl( \bigl\{ |\tilde{\sigma}^1\cdot\tilde{\rho}^1| \leq \sqrt{c_1c_2}+\eps \bigr\}\setminus
\bigl\{ |\tilde{\sigma}^1\cdot\tilde{\rho}^1 -\chi| \leq \delta\bigr\}
\Bigr) \Bigr\ra = 0.
\label{eqabove2}
\end{equation}
\end{proposition}
Together with (\ref{eqabove}), (\ref{eqabove2}) implies that, for any $\delta>0$,
\begin{equation}
\lim_{N\to\infty}\e \Bigl\la I\bigl( 
|\tilde{\sigma}^1\cdot\tilde{\rho}^1 -\chi| \geq \delta
\bigr) \Bigr\ra = 0,
\end{equation}
and this finishes the proof of Theorem \ref{Th1}.

\end{document}